\input epsf

\documentclass[12pt]{article}

\usepackage{latexsym,amsmath,amssymb,amsfonts,amscd,multirow}
\usepackage{epsfig,subfigure,verbatim,epstopdf,graphics}
\usepackage{color,leftidx}
\usepackage{slashbox}
\newtheorem{thm}{Theorem}[section]

\newtheorem{lem}[thm]{Lemma}
\newtheorem{exa}[thm]{Example}

\newcommand{\bx}{\mathbf{x}}

\newcommand{\bE}{\mathbf{E}}

\newfont{\iams}{msbm9}

\newcommand{\commentbis}[1]{}
\newcommand{\be}{\begin{eqnarray}}
\newcommand{\ee}{\end{eqnarray}}
\newcommand{\beno}{\begin{eqnarray*}}
\newcommand{\eeno}{\end{eqnarray*}}
\newcommand{\barr}[1]{\begin{array}{#1}}
\newcommand{\earr}{\end{array}}

\newcommand{\df}{\partial}

\newcommand{\mE}{{\mathcal E}}
\def\Real{\mathbb R}
\newcommand{\beq}{\begin{equation}}
\newcommand{\eeq}{\end{equation}}
\newcommand{\beqa}{\begin{eqnarray}}
\newcommand{\eeqa}{\end{eqnarray}}

\newcommand{\Kx}{{K_x}}
\newcommand{\Kv}{{K_v}}
\newcommand{\bv}{{\bf v}}
\newcommand{\bB}{{\bf B}}
\newcommand{\bJ}{{\bf J}}

\newcommand{\bU}{{\bf U}}

\newcommand{\bn}{{\bf n}}
\newcommand{\Ox}{{\Omega_x}}
\newcommand{\Ov}{{\Omega_v}}

\newcommand{\nh}{{n+1/2}}
\newcommand{\na}{{n+1}}
\newcommand{\mT}{{\mathcal T}}
\newcommand{\mG}{{\mathcal G}}
\newcommand{\mS}{{\mathcal S}}
\newcommand{\mU}{{\mathcal U}}
\newcommand{\mW}{{\mathcal W}}
\newcommand{\mZ}{{\mathcal Z}}
\newcounter{nsez}

\title
{  Energy-conserving discontinuous Galerkin methods for the Vlasov-Amp\`{e}re system}

\author{  Yingda Cheng
\thanks{Department of Mathematics, Michigan State University,
East Lansing, MI 48824 U.S.A.
 {\tt ycheng@math.msu.edu}}
  \and
Andrew J. Christlieb
\thanks{Department of Mathematics, Michigan State University,
East Lansing, MI 48824 U.S.A.
 {\tt christlieb@math.msu.edu}}
   \and
Xinghui Zhong
\thanks{Department of Mathematics, Michigan State University,
East Lansing, MI 48824 U.S.A.
 {\tt zhongxh@math.msu.edu}}
}

\date{\today}

\begin{document}


\maketitle

\begin{abstract}
In this paper, we propose energy-conserving numerical schemes for the Vlasov-Amp\`{e}re (VA) systems. The VA system is a model used to describe the evolution of probability density function of charged particles under self consistent electric field in plasmas. It conserves many physical quantities, including the total energy which is comprised of the kinetic and electric energy. Unlike the total particle number conservation, the total energy conservation is challenging to achieve.  For simulations in longer time ranges, negligence of this fact could cause unphysical results, such as plasma self heating or cooling.
In this paper, we develop the first Eulerian solvers that can preserve fully discrete total energy conservation. The main components of our solvers include explicit or implicit energy-conserving temporal discretizations, an energy-conserving operator splitting for the VA equation and discontinuous Galerkin finite element methods for the spatial discretizations.
We validate our schemes by rigorous derivations and benchmark numerical examples such as Landau damping, two-stream instability and bump-on-tail instability.

\end{abstract}



\section{Introduction}

Plasma  is a state of matter similar to gas in which a certain portion of the particles is ionized. Understanding the complex behavior of plasmas has led to important advances ranging from space physics, fusion energy, to high-power microwave generation and large scale particle  accelerators.
One of the fundamental models in plasma physics is the Vlasov equation, which is a kinetic equation that describes the time  evolution of the probability distribution function of collisionless charged particles with long-range interactions. In plasma physics, those interactions are described by a self-consistent collective field, and can be modeled by the Maxwell's equation or the Poisson's equations in the non-relativistic zero-magnetic field limit, resulting in the  popular Vlasov-Maxwell (VM) or Vlasov-Poisson (VP) systems.

 As for Vlasov solvers, the particle-in-cell (PIC) methods  \cite{Birdsall_book1991, Hockney_book1981} have long been very popular  numerical tools. In PIC methods, the macro-particles are advanced in a Lagrangian framework, while the field equations are solved on a mesh.
 In recent years, there has been growing interest in computing kinetic equations in a deterministic fashion (i.e. the direct computation for the solutions to the Vlasov equations under Eulerian or semi-Lagrangian framework).  Deterministic solvers enjoy the advantage of producing highly accurate results without having any statistical noise.
The main computational challenges for those methods include the high-dimensionality  of the kinetic equation, multiple temporal and spatial scales associated with various physical phenomena,  the conservation of the physical quantities due to the Hamiltonian structure of the systems.
 Deterministic numerical schemes for the VP equations include    semi-Lagrangian methods  \cite{chengknorr_76, Sonnendrucker_1999, qiu2010conservative,rossmanith2011positivity, qiu2011positivity,PhysRevA.39.6356,PhysRevLett.63.2361,PhysRevA.43.4452}, the
WENO method coupled with Fourier collocation \cite{Guo_landau_2001}, finite volume  methods  \cite{Boris_1976,Fijalkow_1999, Filbet_PFC_2001}, Fourier-Fourier spectral methods \cite{Klimas_1987, Klimas_1994}, continuous finite element methods \cite{Zaki_1988_1, Zaki_1988_2}, Runge-Kutta discontinuous Galerkin (DG) methods \cite{heath2012discontinuous, Heath_thesis, Ayuso2009, Ayuso2010, cheng_vp}, among many others.
As for VM systems,  PIC methods \cite{brackbill1982implicit, Jacobs_Hesthaven_2006, Jacobs_Hesthaven_2009, Markidis20117037},   semi-Lagrangian methods \cite{califano1998ksw, mangeney2002nsi, califano1965ikp, califano2001ffm}, spectral methods \cite{eliasson2003numerical}, finite difference methods \cite{umeda2009two}, and RKDG methods \cite{cheng_vm} have been developed for many applications.

In this paper, we focus on the Vlasov-Amp\`{e}re (VA) systems. The VA system could be viewed as the zero-magnetic limit of the VM system, and is  equivalent to the VP system if charge continuity is satisfied, when there is no external applied voltage. Unlike the VP system, the electric field is not solved from the Poisson equation, but is advanced in time through the current generated by the moving charges. Therefore the algorithms designed for the VA system has important implications for the design of VM solvers.  A PIC algorithm \cite{chen2011energy},  a finite difference method \cite{Horne2001182}, a semi-Lagrangian method \cite{crouseilles_va}, and a finite volume scheme \cite{elkina2006new} have been proposed to solve the VA system.
In this paper, we design Eulerian solvers to treat the VA equations. The   proposed  schemes have several important features: it conserves the total particle number and energy of the system on the fully discrete level;    it has a systematic way to incorporate explicit or implicit  time stepping depending on the stiffness of the equations; and it could be designed for implementations on unstructured grids for complex geometries in the physical space.

One particular important focus of this paper is energy conservation. For most methods in the  Vlasov literature, the conservation of the total particle number is maintained, but  the conservation of total energy is not addressed, rather it was left to the accuracy of the scheme. For simulations in longer time ranges, the spurious energy created or annihilated by numerical methods could build up and lead to unphysical results, such as plasma self heating or cooling \cite{cohen1989performance}. This issue will be more prominent if we use under-resolved mesh or large time steps. In \cite{PhysRevA.39.6356,PhysRevLett.63.2361,PhysRevA.43.4452}, the authors developed energy conserving convective schemes for plasma simulations.
Recently, several PIC methods have been proposed to conserve the total energy. In \cite{chen2011energy}, PIC for VA equations was developed; it is fully implicit, energy and charge conserving. In \cite{Markidis20117037}, PIC for VM system was developed, in which Maxwell equations is solved on Yee's lattice \cite{yee1966numerical} and implicit midpoint method  is used as the time integrator. In \cite{Filbet_compare_2003, ayuso2012high},  finite difference and DG methods were proposed to conserve the total energy of VP systems on the semi-discrete level.  In this paper, we develop the first Eulerian solvers that can conserve the total energy on the fully discrete level, and they incorporate explicit and implicit time steppings and are suitable for implementation on unstructured grids for complex geometries.

The rest of the paper is  organized as follows: in Section \ref{sec:equation}, we describe the models under consideration.
In Section \ref{sec:method}, the schemes and their  properties are discussed.   Section \ref{sec:numerical} is devoted to numerical results. We conclude with a few remarks in Section \ref{sec:conclusion}.


\section{Basic Equations}
\label{sec:equation}

In this section, we will introduce the basic models. We consider the evolution of electron probability distribution function in the presence of a uniform background of ions.
Under the scaling of the characteristic time by the inverse of the plasma frequency $\omega_p^{-1}$ and length scaled by the Debye length $\lambda_D$, and characteristic electric and magnetic field as $\bar{E}=\bar{B}=-m c \omega_p/e$, the VM system is formulated as follows
 \begin{eqnarray}
&&\partial_t f  +  \bv \cdot \nabla_{\bx} f  + (\bE + \bv \times \bB) \cdot \nabla_\bv f = 0~, \notag\\
&&\frac{\df \bE}{\df t} = \nabla_\bx \times \bB - \bJ, \qquad
\frac{\df \bB}{\df t} = -  \nabla_\bx \times \bE~,  \quad\label{modelvm} \\
&&\nabla_\bx \cdot \bE =  \rho-\rho_i, \qquad
\nabla_\bx \cdot \bB = 0~, \quad\notag
\end{eqnarray}
where the density and current density are defined as
\begin{equation*}
\rho(\bx, t)= \int_\Ov f(\bx, \bv, t)d\bv,\qquad \bJ(\bx, t)= \int_\Ov f(\bx, \bv, t)\bv d\bv.
\end{equation*}
and $\rho_i$ is the ion density.
In this model, $f=f(\bx,\bv,t)$ is  the probability distribution function ($pdf$) for  finding an electron (at position $x$ with velocity $v$ at time $t$) with a uniform background of fixed ions under a self-consistent electrostatic field.   Here the domain $\Omega= \Omega_x \times \Ov$, where $\Omega_x$ can be either a finite domain or $\Real^n$ and $\Ov=\Real^n$. The boundary conditions for the above systems are summarized as follows:  $f\rightarrow0$ as $|x|\rightarrow \infty$ or $|v|\rightarrow \infty$. If $\Omega_x$ is finite, then we can  impose either inflow boundary conditions with $f=f^{in}$ on $\Gamma_I=\{(x,v)| v \cdot \nu_x < 0\}$, where $\nu_x$ is the outward normal vector,  or more simply impose  periodic boundary conditions. For simplicity of discussion, in this paper, we will always assume periodicity in  $x$.  In practice the velocity domain needs to be truncated, so that $\Ov$ is finite. Discussions related to the domains could be found for example in \cite{cheng_vp, cheng_vm}. In particular, in all the discussions below, we will assume that $\Ov$ is taken large enough, so that the numerical solution $f_h \approx 0$ at $\partial \,\Ov$. We could do this by enlarging the velocity domain enough and some related discussions could be found in \cite{cheng_vm}.

In the zero-magnetic limit, the VM system becomes
\beqa
\label{vap}
&\partial_t f  +  \bv \cdot \nabla_{\bx} f  + \bE \cdot \nabla_\bv f = 0~, \\
&\frac{\df \bE}{\df t} = - \bJ, \qquad
\nabla_\bx \cdot \bE =  \rho-\rho_i, \quad\notag
\eeqa
This leads to either the Vlasov-Amp\`{e}re (VA) system
\beqa
\label{va}
&\partial_t f  +  \bv \cdot \nabla_{\bx} f  + \bE \cdot \nabla_\bv f = 0~,  \\
&\frac{\df \bE}{\df t} = - \bJ, \quad\notag
\eeqa
or the VP system
\beqa
\label{vp}
&\partial_t f  +  \bv \cdot \nabla_{\bx} f  + \bE \cdot \nabla_\bv f = 0~, \\\
&\nabla_\bx \cdot \bE =  \rho-\rho_i, \quad\notag
\eeqa
 In the absence of external fields, the VA and VP systems are equivalent when the charge continuity equation
$$\rho_t+\nabla_\bx \cdot \bJ=0$$
is satisfied. The initial condition of the electric field $\bE$ for the VA system can be provided by solving the Poisson equation.
It is well-known that the VA and VP systems conserves the total particle number $\int_\Omega f \, d\bx d\bv$,   and  the total energy $$TE=\frac{1}{2} \int _\Ox \int_\Ov f |\bv|^2 d\bv d\bx + \frac{1}{2}\int_\Ox |\bE|^2 d\bx, $$
which is comprised of the kinetic and electric energy. Moreover, any functional of the form $\int _\Ox \int_\Ov G(f)\, d\bv d\bx $ is a constant of motion. In particular, this includes the $k$-th order invariant $I_k = \int _\Ox \int_\Ov f^k \,d\bv d\bx $ and the entropy $S= -\int _\Ox \int_\Ov f \ln(f) \, d\bv d\bx $.  Sometimes the functional $I_2$ is also called the enstrophy, and all of these invariants are called Casimir invariants (see, e.g., \cite{morrison_98}). The total momentum of the system can be defined by $\int _\Ox \int_\Ov f \bv \,d\bv d\bx$ and this is conserved for many well known examples, such as Landau damping and two-stream instability.

\section{Numerical Methods}
\label{sec:method}

In this section, we will describe the numerical methods for the VA system \eqref{va} and discuss their properties. The main components of the proposed schemes include energy-conserving temporal and spatial discretizations.

\subsection{Temporal discretizations}
In this subsection, we will describe several versions of temporal discretizations, while leaving the  $(\bx, \bv)$ variables continuous. We will first establish second-order explicit and implicit schemes. Then to treat the implicit methods more efficiently without inverting   nonlinear systems in the whole $(\bx, \bv)$ space, we propose a split-in-time algorithm. In this algorithm, the VA system is splitted into two equations, each of which conserves the total energy, and could be solved in reduced dimensions. Finally, we  discuss how to generalize the methods beyond second order.

\subsubsection{Second order schemes}
\label{sec:secondorder}
An explicit second order scheme can be designed as follows
\begin{subequations}
\label{scheme1}
\begin{align}
&\frac{f^{n+1/2}-f^n}{\triangle t/2}   +\bv \cdot \nabla_{\bx} f^n  +  \bE^n \cdot \nabla_\bv f^n = 0~, \label{scheme1:1}\\
&\frac{\bE^{n+1}-\bE^n}{\triangle t}=-\bJ^{n+1/2},  \quad \textrm{where}\, \,\bJ^{n+1/2}= \int_\Ov f^{n+1/2} \bv d\bv \label{scheme1:2}\\
&\frac{f^{n+1}-f^n}{\triangle t}  + \bv \cdot \nabla_{\bx} f^{n+1/2}   +  \frac{1}{2}(\bE^n+\bE^{n+1})  \cdot \nabla_\bv f^{n+1/2} = 0~. \label{scheme1:3}
\end{align}
\end{subequations}
The key idea is the careful coupling of the Vlasov and Amp\`{e}re solvers. We will verify the conservation property in Theorem \ref{thm:energy}. We denote the scheme above to be {\bf Scheme-1}$(\triangle t)$, namely, this means $(f^\na, \bE^\na)= \textnormal{\bf Scheme-1}(\triangle t) (f^n, \bE^n)$.

A fully implicit scheme based on the implicit midpoint method can also preserve the total energy.
\begin{subequations}
\label{scheme2}
\begin{align}
&\frac{f^{n+1}-f^n}{\triangle t}   + \bv \cdot \nabla_{\bx} \frac{f^n+f^{n+1}}{2} +  \frac{1}{2}(\bE^n+\bE^{n+1})  \cdot \nabla_\bv \frac{f^n+f^{n+1}}{2} = 0~,\label{scheme2:1}\\
&\frac{\bE^{n+1}-\bE^n}{\triangle t}=- \frac{1}{2} (\bJ^{n}+\bJ^{n+1}). \label{scheme2:2}
\end{align}
\end{subequations}

The scheme above can remove the CFL restriction of the explicit scheme, but it involves nonlinearly coupled computations of $\bE$ and $f$. We denote the scheme above by $\textnormal{\bf Scheme-2}(\triangle t)$.

 On the other hand, we could   construct by Strang splitting another implicit scheme in which the solver for $\bE$ and $f$ is decoupled, and we denote it to be $\textnormal{\bf Scheme-3}(\triangle t)$.
\begin{subequations}
\label{scheme3}
\begin{align}
&\frac{\bE^{n+1/2}-\bE^n}{\triangle t/2}=-\bJ^{n},  \label{scheme3:1}\\
&\frac{f^{n+1}-f^n}{\triangle t}  + \bv \cdot \nabla_{\bx} \frac{f^n+f^{n+1}}{2}  +  \bE^\nh  \cdot \nabla_\bv \frac{ f^n+f^\na}{2} = 0~,\label{scheme3:2}\\
&\frac{\bE^{n+1}-\bE^\nh}{\triangle t/2}=-\bJ^{n+1}.\label{scheme3:3}
\end{align}
\end{subequations}

Through simple Taylor expansions, we can verify that all three schemes above are  second order accurate in time. The implicit schemes \eqref{scheme2}, \eqref{scheme3}  are also symmetric in time (time reversible). In the next theorem, we will prove energy conservation for the methods above.

\begin{thm}[{Total energy conservation}]
\label{thm:energy}
With the boundary conditions described in Section 2, the schemes above preserve the discrete total energy $TE_n=TE_{n+1}$, where
$$2\,  (TE_n)=\int_\Ox \int_\Ov f^n |\bv|^2 d\bv d\bx +\int_\Ox |\bE^n|^2 d\bx$$
in {\bf Scheme-1}$(\triangle t)$ and  {\bf Scheme-2}$(\triangle t)$, and
\begin{eqnarray*}
2\,  (TE_n)&=&\int_\Ox \int_\Ov f^n |\bv|^2 d\bv d\bx +\int_\Ox \bE^{n+1/2} \cdot \bE^{n-1/2} d\bx\\
&=&\int_\Ox \int_\Ov f^n |\bv|^2 d\bv d\bx +\int_\Ox |\bE^n|^2 d\bx-\frac{\triangle t^2}{4} \int_\Ox |\bJ^n|^2 d\bx
\end{eqnarray*}
in {\bf Scheme-3}$(\triangle t)$.
\end{thm}

\emph{Proof.} As for  $\textnormal{\bf Scheme-1}(\triangle t)$,
\begin{eqnarray*}
&&\int_\Ox \int_\Ov \frac{f^{n+1}-f^n}{\triangle t} |\bv|^2 d\bv d\bx \\
&&=-\int_\Ox \int_\Ov \bv \cdot \nabla_{\bx} f^{n+1/2}  |\bv|^2 d\bv d\bx   -  \int_\Ox \int_\Ov   \frac{1}{2}(\bE^n+\bE^{n+1})  \cdot \nabla_\bv f^{n+1/2}  |\bv|^2 d\bv d\bx  \\
&& =-\int_\Ov  |\bv|^2 \bv \cdot (\int_\Ox  \nabla_{\bx} f^{n+1/2}  d\bx)\, d\bv   +  \int_\Ox \int_\Ov  (\bE^n+\bE^{n+1})\cdot \left( f^{n+1/2}  \bv\right) d\bv d\bx \\
&&=\int_\Ox   (\bE^n+\bE^{n+1})\cdot \,\bJ^{n+1/2}   d\bx
\end{eqnarray*}
and
\begin{eqnarray*}
&&\int_\Ox  \frac{\bE^{n+1}-\bE^n}{\triangle t} \cdot (\bE^{n+1}+\bE^n)d\bx=-\int_\Ox   (\bE^n+\bE^{n+1}) \cdot\bJ^{n+1/2}   d\bx.
\end{eqnarray*}
Therefore,
$$\int_\Ox \int_\Ov f^n |\bv|^2 d\bv d\bx +\int_\Ox |\bE^n|^2 d\bx=\int_\Ox \int_\Ov f^{n+1} |\bv|^2 d\bv d\bx +\int_\Ox |\bE^{n+1}|^2 d\bx.$$
The other two proof are similar and are omitted.
 $\Box$
 
 From this theorem, we can see that  {\bf Scheme-1} and {\bf Scheme-2} exactly preserve the total energy, while {\bf Scheme-3} achieves near conservation of the total energy. The numerical energy from {\bf Scheme-3} is a second order modified version of the original total energy. This is natural due to the second order accuracy of the scheme.

\subsubsection{Split in time}

The implicit schemes   require inverting the problem in $(\bx, \bv)$ space, which is costly for high-dimensional applications. In this subsection, we propose a splitting framework for the VA equations so that the resulting equations could be computed in reduced dimensions.
Using splitting schemes to treat the Vlasov equation is not new. In fact, the very popular semi-Lagrangian methods for Vlasov simulations are based on the dimensional splitting of the Vlasov equation \cite{chengknorr_76}. In those methods, the Vlasov equation is splitted into several one-dimensional equations, which becomes advection equations that the semi-Lagrangian methods could easily handle.
Here we propose to split not just the Vlasov equation  but the whole VA system together, and each splitted equation can still maintain energy conservation.

For the model VA equation \eqref{va}, we propose to perform the operator splitting as follows:

\[
  \textrm{ (a)} \left\{
  \begin{array}{l}
\partial_t f  +  \bv \cdot \nabla_{\bx} f   = 0~, \\
\partial_t \bE = 0,
  \end{array} \right.
  \textrm{(b)} \left\{
  \begin{array}{l}
\partial_t f  +   \bE  \cdot \nabla_\bv f   = 0~, \\
\partial_t \bE = -\bJ,
  \end{array} \right.
\]
One of the main feature of this splitting is that each of the two equations is energy-conserving,
$$\frac{d}{dt} (\int_\Ox \int_\Ov f |\bv|^2 d\bv d\bx+ \int_\Ox |\bE|^2 d\bx ) = 0.$$ In particular,
\[
  \textrm{ (a)} \left\{
  \begin{array}{l}
\displaystyle\frac{d}{dt} \int_\Ox \int_\Ov f |\bv|^2 d\bv d\bx = 0~, \\[5mm]
\displaystyle\frac{d}{dt}  \int_\Ox |\bE|^2 d\bx = 0,
  \end{array} \right.
  \textrm{(b)} \frac{d}{dt} (\int_\Ox \int_\Ov f |\bv|^2 d\bv d\bx+ \int_\Ox |\bE|^2 d\bx ) = 0~,
\]
We can see that equation (a) contains the free streaming operator. Equation (b) contains the interchange of kinetic and electric energy. Therefore, we only need to design energy-conserving temporal discretizations for each equations and then carefully couple the two solvers together to achieve the conservation for the whole system, and the desired accuracy.

As for equation (a), we can use any implicit or explicit Runge-Kutta methods to solve it, and they all conserve the kinetic energy. To see this, consider the forward Euler
$$\frac{f^{n+1}-f^n}{\triangle t} +  \bv \cdot \nabla_{\bx} f^n=0,$$
or  backward Euler method
$$\frac{f^{n+1}-f^n}{\triangle t} +  \bv \cdot \nabla_{\bx} f^\na=0,$$
A simple check yields
$\int_\Ox \int_\Ov f^n |\bv|^2 d\bv d\bx =\int_\Ox \int_\Ov f^\na |\bv|^2   d\bv d\bx $. (Note that here we have abused the notation, and use superscript $n$, $n+1$ to denote the sub steps in computing equation (a), not the whole time step to compute the VA system). Therefore, we could pick a suitable Runge-Kutta method with desired order and property for this step. To be second order, one could use the implicit midpoint method,
\beq
\label{schemea}
\frac{f^{n+1}-f^n}{\triangle t} +  \bv \cdot \nabla_{\bx} \frac{f^n+f^\na}{2}=0.
\eeq

Equation (b) contains the main coupling effect of the Vlasov and Amp\`{e}re equations, and has to be computed carefully to ensure balance of kinetic and electric energy. We could use the methods studied in  Section \ref{sec:secondorder} to compute this equation. (We only need to include the corresponding terms as those appeared in equation (b)).  The resulting scheme will naturally preserve a discrete form of  the sum of kinetic and electric energy.

Scheme \eqref{scheme2},  the implicit midpoint method will reduce to
\begin{subequations}
\label{scheme2s}
\begin{align}
&\frac{f^{n+1}-f^n}{\triangle t}    +  \frac{1}{2}(\bE^n+\bE^{n+1})  \cdot \nabla_\bv \frac{f^n+f^{n+1}}{2} = 0~,\label{scheme2s:1}\\
&\frac{\bE^{n+1}-\bE^n}{\triangle t}=- \frac{1}{2} (\bJ^{n}+\bJ^{n+1}) \label{scheme2s:2}
\end{align}
\end{subequations}

Finally, suppose we use $\textnormal{\bf Scheme-a}(\triangle t)$ to denote  second order schemes for equation (a), and $\textnormal{\bf Scheme-b}(\triangle t)$ to denote second order schemes for equation (b), then we can show that by Strang splitting 
$\textnormal{\bf Scheme-4}(\triangle t)=\textnormal{\bf Scheme-a}(\triangle t/2) \textnormal{\bf Scheme-b}(\triangle t) \textnormal{\bf Scheme-a}(\triangle t/2)$, the method is second order for the original VA system.

\begin{thm}[{Total energy conservation for the splitted methods}]
\label{thm:energy2}
With the boundary conditions described in Section 2, the scheme with {\bf Scheme-a} being \eqref{schemea} and {\bf Scheme-b} being \eqref{scheme2s}, and $\textnormal{\bf Scheme-4}(\triangle t):=\textnormal{\bf Scheme-a}(\triangle t/2) \textnormal{\bf Scheme-b}(\triangle t) \textnormal{\bf Scheme-a}(\triangle t/2)$
 preserves the discrete total energy $TE_n=TE_{n+1}$, where
$$2\,  (TE_n)=\int_\Ox \int_\Ov f^n |\bv|^2 d\bv d\bx +\int_\Ox |\bE^n|^2 d\bx.$$
\end{thm}
The proof is straightforward by the discussion in this subsection and is omitted.

\subsubsection{Generalizations to higher order}
\label{sec:high}
We could generalize the   second order schemes to higher order. High order symplectic methods were constructed in \cite{yoshida1990construction, forest1990fourth}, and  how to generalize second order time reversible schemes into fourth order time symmetric schemes has been demonstrated in \cite{de1992easily, sanz1991order}.

Using the idea of  \cite{de1992easily}, we let $\beta_1, \beta_2, \beta_3$ satisfy
$$\beta_1+\beta_2+\beta_3=1,  \qquad \beta^3_1+\beta^3_2+\beta^3_3=1, \qquad \beta_1=\beta_3,$$
from which we get $\beta_1=\beta_3=(2+2^{1/3}+2^{-1/3})/3\approx1.3512$, $\beta_2=1-2\beta_1\approx-1.7024$.

Then we let
$$\textnormal{\bf Scheme-2F}(\triangle t)=\textnormal{\bf Scheme-2}(\beta_1 \triangle t) \textnormal{\bf Scheme-2}(\beta_2 \triangle t) \textnormal{\bf Scheme-2}(\beta_3 \triangle t);$$
$$\textnormal{\bf Scheme-3F}(\triangle t)=\textnormal{\bf Scheme-3}(\beta_1 \triangle t) \textnormal{\bf Scheme-3}(\beta_2 \triangle t) \textnormal{\bf Scheme-3}(\beta_3 \triangle t);$$
$$\textnormal{\bf Scheme-4F}(\triangle t)=\textnormal{\bf Scheme-4}(\beta_1 \triangle t) \textnormal{\bf Scheme-4}(\beta_2 \triangle t) \textnormal{\bf Scheme-4}(\beta_3 \triangle t).$$
 We can verify that $\textnormal{\bf Scheme-2F}(\triangle t), \textnormal{\bf Scheme-3F}(\triangle t),
 \textnormal{\bf Scheme-4F}(\triangle t)$ are all fourth order. For {\bf Scheme-1}$(\triangle t)$, because it is not symmetric in time, this procedure will not be able to raise the method to fourth order accuracy.

\begin{thm}[{Total energy conservation for the fourth order methods}]
\label{thm:energy2}
With the boundary conditions described in Section 2, 
{\bf Scheme-2F}$(\triangle t)$ and {\bf Scheme-4F}$(\triangle t)$
 preserves the discrete total energy $TE_n=TE_{n+1}$, where
$$2\,  (TE_n)=\int_\Ox \int_\Ov f^n |\bv|^2 d\bv d\bx +\int_\Ox |\bE^n|^2 d\bx.$$
\end{thm}
The proof is straightforward by the properties of the second order methods {\bf Scheme-2} and {\bf Scheme-4} and is omitted.

 The theorem above shows that $\textnormal {\bf Scheme-2F}(\triangle t)$ and $\textnormal{\bf Scheme-4F}(\triangle t)$ exactly preserve
 the total energy.  On the other hand, it is  challenging to obtain the explicit form of the modified total energy for $\textnormal{\bf Scheme-3F}(\triangle t)$. In Section \ref{sec:numerical}, we demonstrate numerically that it can achieve near conservation of the total energy.

Finally, we remark that since $\beta_2<0$,  special care needs to be taken   for those negative time steps, e.g. the flux discussed in the next subsection needs to be reversed, i.e. upwind flux needs to be replaced by downwind flux.

\subsection{Fully discrete methods}
In this section, we will discuss the spatial discretizations and formulate the fully discrete schemes. In particular, we consider two approaches: one being the unsplit schemes, the other being the splitted implicit schemes.

In this paper, we choose to use discontinuous Galerkin (DG) methods to discretize the $(\bx, \bv)$ variable due to their excellent conservation properties.
The DG method \cite{Cockburn_2000_history, Cockburn_2001_RK_DG} is a class
of finite element methods using discontinuous piecewise
polynomial space for the numerical solution and the test functions, and they have excellent conservation properties.  DG methods have been designed to solve VP \cite{heath2012discontinuous, Heath_thesis, Ayuso2009, Ayuso2010, cheng_vp, rossmanith2011positivity, qiu2011positivity} and VM \cite{cheng_vm} systems. In particular, semi-discrete total energy conservation have been established in \cite{ Ayuso2009,ayuso2012high} for  VP and in \cite{cheng_vm} for VM systems. In the discussions below, we will prove fully discrete conservation properties for our proposed methods.

\subsubsection{Notation}

Let $\mT_h^x=\{\Kx\}$ and $\mT_h^v=\{\Kv\}$ be  partitions of $\Ox$ and $\Ov$, respectively,  with $\Kx$ and $\Kv$ being  Cartesian elements or simplices;  then $\mT_h=\{K: K=\Kx\times\Kv, \forall \Kx\in\mT_h^x, \forall \Kv\in\mT_h^v\}$ defines a partition of $\Omega$. Let $\mE_x$ be the set of the edges of $\mT_h^x$ and  $\mE_v$ be the set of the edges of $\mT_h^v$;  then the edges of $\mT_h$ will be $\mE=\{\Kx\times e_v: \forall\Kx\in\mT_h^x, \forall e_v\in\mE_v\}\cup \{e_x\times\Kv: \forall e_x\in\mE_x, \forall \Kv\in\mT_h^v\}$. Here we take into account the periodic boundary condition in the $\bx$-direction when defining $\mE_x$ and $\mE$. Furthermore, $\mE_v=\mE_v^i\cup\mE_v^b$ with $\mE_v^i$ and $\mE_v^b$ being the set of interior and boundary edges of $\mT_h^v$, respectively.

We will make use of  the following discrete spaces
\begin{subequations}
\begin{align}
\mG_h^{k}&=\left\{g\in L^2(\Omega): g|_{K=\Kx\times\Kv}\in P^k(\Kx\times\Kv), \forall \Kx\in\mT_h^x, \forall \Kv\in\mT_h^v \right\}~,\label{eq:sp:f}\\
\mS_h^{k}&=\left\{g\in L^2(\Omega): g|_{K=\Kx\times\Kv}\in P^k(\Kx)\times P^k(\Kv), \forall \Kx\in\mT_h^x, \forall \Kv\in\mT_h^v \right\}~,\\
\mU_h^k&=\left\{\bU\in [L^2(\Omega_x)]^{d_x}: \bU|_\Kx\in [P^k(\Kx)]^{d_x}, \forall \Kx\in\mT_h^x \right\}~,\\
\mW_h^k&=\left\{w\in L^2(\Omega_x): w|_\Kx\in P^k(\Kx), \forall \Kx\in\mT_h^x \right\}~,\\
\mZ_h^k&=\left\{z\in L^2(\Omega_v): w|_\Kv\in P^k(\Kv), \forall \Kv\in\mT_h^v \right\}~,
\end{align}
\end{subequations}
where $P^k(D)$ denotes the set of polynomials of  total degree at most $k$ on $D$.
The discussion about those spaces for Vlasov equations can be found in \cite{cheng_vp, cheng_vm}.

For piecewise  functions defined with respect to $\mT_h^x$ or $\mT_h^v$, we further introduce the jumps and averages as follows. For any edge $e=\{K_x^+\cap K_x^-\}\in\mE_x$, with $\bn_x^\pm$ as the outward unit normal to $\partial K_x^\pm$,
$g^\pm=g|_{K_x^\pm}$, and $\bU^\pm=\bU|_{K_x^\pm}$, the jumps across $e$ are defined  as
\begin{equation*}
[g]_.={g^+}{\bn_.^+}+{g^-}{\bn_.^-},\qquad [\bU]_.={\bU^+}\cdot{\bn_.^+}+{\bU^-}\cdot{\bn_.^-}
\end{equation*}
and the averages are
\begin{equation*}
\{g\}_{.}=\frac{1}{2}({g^+}+{g^-}),\qquad \{\bU\}_.=\frac{1}{2}({\bU^+}+{\bU^-}),
\end{equation*}
where $.$ are used to denote $\bx$ or $\bv$.

\subsubsection{Unsplit schemes and their properties}

In this subsection, we will describe the DG methods for the unsplit schemes {\bf Scheme-1}, {\bf Scheme-2}, {\bf Scheme-3}, {\bf Scheme-2F}, {\bf Scheme-3F}, and discuss their properties.  For example, the scheme with {\bf Scheme-1}$(\triangle t)$  is formulated as follows: we look for $ f_h^{n+1/2},f_h^{n+1}\in\mG_h^k$,
such that for any $g\in\mG_h^k$,
\begin{subequations}
\label{dscheme1}
\begin{align}
&\int_K \frac{f_h^{n+1/2}-f_h^n}{\triangle t/2}  g d\bx d\bv
- \int_K f_h^n\bv\cdot\nabla_\bx g d\bx d\bv
- \int_K f_h^n \bE_h^n \cdot\nabla_\bv g d\bx d\bv\notag\\[3mm]
&\quad+ \int_{\Kv}\int_{\df\Kx} \widehat{f_h^n \bv\cdot \bn_x} g ds_x d\bv + \int_{\Kx} \int_{\df\Kv} \widehat{(f_h^n \bE_h^n \cdot \bn_\bv)} g ds_v d\bx=0~,\label{dscheme1:1}\\[3mm]
& \frac{\bE_h^{n+1}-\bE_h^n}{\triangle t}=-  \bJ_h^{n+1/2}\,  \quad \textrm{where}\, \,\bJ_h^{n+1/2}= \int_\Ov f_h^{n+1/2} \bv d\bv~,\label{dscheme1:2}\\[3mm]
& \int_K \frac{f_h^{n+1}-f_h^n}{\triangle t}  g d\bx d\bv
- \int_K f_h^{n+1/2}\bv\cdot\nabla_\bx g d\bx d\bv
- \frac{1}{2} \int_K f_h^{n+1/2} (\bE_h^n+\bE_h^{n+1}) \cdot\nabla_\bv g d\bx d\bv\notag\\[3mm]
&\quad+ \int_{\Kv}\int_{\df\Kx} \widehat{f_h^{n+1/2} \bv\cdot \bn_x} g ds_x d\bv + \frac{1}{2}\int_{\Kx} \int_{\df\Kv} \widehat{(f_h^{n+1/2} (\bE_h^n+\bE_h^{n+1}) \cdot \bn_\bv)} g ds_v d\bx=0~,\label{dscheme1:3}
\end{align}
\end{subequations}
Here $\bn_x$ and $\bn_v$ are outward unit normals of $\df\Kx$ and $\df\Kv$, respectively.  Since $\bJ_h^{n+1/2} \in \mU_h^k$, the space that $\bE_h^n$ lies in are totally determined by $\bE_h^0$, i.e. the initial electric field. All  ``hat''  functions are numerical fluxes that are determined by  either upwinding, i.e.,
\begin{subequations}
\begin{align}
\widehat{f_h^n \bv\cdot \bn_x}:&=\widetilde{f_h^n \bv}\cdot \bn_x=\left(\{f_h^n\bv\}_x+\frac{|\bv\cdot\bn_x|}{2}[f_h^n]_x\right)\cdot\bn_x~,\label{eq:flux:1}\\
\widehat{f_h^n \bE_h^n \cdot \bn_\bv}:&=\widetilde{f_h^n \bE_h^n }\cdot\bn_\bv=\left(\{f_h^n \bE_h^n\}_\bv+\frac{|\bE_h^n\cdot\bn_\bv|}{2}[f_h^n]_\bv\right)\cdot\bn_\bv~,
\end{align}
\end{subequations}
or central flux
\begin{subequations}
\begin{align}
\widehat{f_h^n \bv\cdot \bn_x}:&=\widetilde{f_h^n \bv}\cdot \bn_x= \{f_h^n\bv\}_x \cdot\bn_x~,\label{eq:flux:1}\\
\widehat{f_h^n \bE_h^n \cdot \bn_\bv}:&=\widetilde{f_h^n \bE_h^n }\cdot\bn_\bv= \{f_h^n \bE_h^n\}_\bv \cdot\bn_\bv~,
\end{align}
\end{subequations}
The upwind and central fluxes in \eqref{dscheme1:3} are defined similarly. It is well known that the upwind flux is more dissipative and the central flux is more dispersive. With the central flux, DG methods for the linear transport equation has sub-optimal order for odd degree polynomials. A numerical comparison of central and upwind fluxes for the VA system is  shown in Section \ref{sec:numerical}. 

The schemes with {\bf Scheme-2}$(\triangle t)$, {\bf Scheme-3}$(\triangle t)$, {\bf Scheme-2F}$(\triangle t)$, {\bf Scheme-3F}$(\triangle t)$ can be formulated similarly, i.e. to use DG discretization to approximate the derivatives of $f$ in $\bx, \bv$. To save space, we do not formulate them here. 

Below we discuss the conservation properties of the fully discrete methods.

\begin{thm}[{Total particle number conservation}]
The scheme \eqref{dscheme1} preserves the total particle number of the system, i.e.
$$\int_\Ox \int_\Ov f_h^\na d\bv d\bx =\int_\Ox \int_\Ov f_h^n  d\bv d\bx .$$
This also holds for DG methods with time integrators {\bf Scheme-2}$(\triangle t)$, {\bf Scheme-3}$(\triangle t)$, {\bf Scheme-2F}$(\triangle t)$, {\bf Scheme-3F}$(\triangle t)$ .
\end{thm}

\emph{Proof.}  Let $g=1$ in \eqref{dscheme1:3}, and sum over all element $K$.  The proof for {\bf Scheme-2}$(\triangle t)$, {\bf Scheme-3}$(\triangle t)$, {\bf Scheme-2F}$(\triangle t)$, {\bf Scheme-3F}$(\triangle t)$  is similar thus omitted. $\Box$

\begin{thm}[{Total energy conservation}]
If $k \geq 2$, the scheme \eqref{dscheme1} preserves the discrete total energy $TE_n=TE_{n+1}$, where
$$2 \,(TE_n)=\int_\Ox \int_\Ov f_h^n |\bv|^2 d\bv d\bx +\int_\Ox |\bE_h^n|^2 d\bx.$$
This also holds for DG methods with time integrators {\bf Scheme-2}$(\triangle t)$, {\bf Scheme-2F}$(\triangle t)$,
and
\begin{eqnarray*}
2\,  (TE_n)&=&\int_\Ox \int_\Ov f_h^n |\bv|^2 d\bv d\bx +\int_\Ox \bE_h^{n+1/2} \cdot \bE_h^{n-1/2} d\bx\\
&=&\int_\Ox \int_\Ov f_h^n |\bv|^2 d\bv d\bx +\int_\Ox |\bE_h^n|^2 d\bx-\frac{\triangle t^2}{4} \int_\Ox |\bJ_h^n|^2 d\bx
\end{eqnarray*}
where $\bJ_h^n= \int_\Ov f_h^{n} \bv d\bv$, is preserved for DG methods with time integrator
{\bf Scheme-3}$(\triangle t)$.
\end{thm}

\emph{Proof.}  Since $|\bv|^2 \in\mG_h^k$, we can let $g=|\bv|^2$ in  \eqref{dscheme1:3}, and sum over all element $K$. Because $g$ is continuous, $\nabla_\bx g=0$, $\nabla_\bv g= 2 \bv$,
\begin{eqnarray*}
\int_\Ox \int_\Ov \frac{f_h^{n+1}-f_h^n}{\triangle t} |\bv|^2 d\bv d\bx
 =  \int_\Ox \int_\Ov  (\bE_h^n+\bE_h^{n+1})\cdot \left( f_h^{n+1/2}  \bv\right) d\bv d\bx
=\int_\Ox   (\bE_h^n+\bE_h^{n+1}) \cdot \bJ_h^{n+1/2}   d\bx
\end{eqnarray*}

and
\begin{eqnarray*}
&&\int_\Ox  \frac{\bE_h^{n+1}-\bE_h^n}{\triangle t} \cdot (\bE_h^{n+1}+\bE_h^n)d\bx=-\int_\Ox   (\bE_h^n+\bE_h^{n+1}) \cdot \bJ_h^{n+1/2}   d\bx
\end{eqnarray*}
Therefore,
$$\int_\Ox \int_\Ov f_h^n |\bv|^2 d\bv d\bx +\int_\Ox |\bE_h^n|^2 d\bx=\int_\Ox \int_\Ov f_h^{n+1} |\bv|^2 d\bv d\bx +\int_\Ox |\bE_h^{n+1}|^2 d\bx,$$
and we are done.

The proof for  {\bf Scheme-2}$(\triangle t)$, {\bf Scheme-2F}$(\triangle t)$ and $ \textnormal{\bf Scheme-3}(\triangle t)$ is similar and thus is omitted.
 $\Box$

Remark: Similar to the discussion in Section \ref{sec:high}, it is  hard to obtain the explicit form of the modified total energy  for $\textnormal{\bf Scheme-3F}(\triangle t)$, and we choose to demonstrate the near conservation of total energy numerically in Section \ref{sec:numerical}.

\begin{thm}[{Charge conservation}]
The scheme \eqref{dscheme1} with central flux satisfies charge continuity. In particular, if the initial electric filed satisfies
\be
\label{e0}
-\int_\Kx \bE_h^0 \cdot \nabla_\bx w d\bx +\int_{\df\Kx} \widehat{\bE_h^0 \cdot \bn_x} w ds_x=\int_\Kx (\rho_h^0-\rho_i) w d\bx,
\ee
for any $w \in \mW_h^k$, then
$$-\int_\Kx \bE_h^n \cdot \nabla_\bx w d\bx +\int_{\df\Kx} \widehat{\bE_h^n \cdot \bn_x} w ds_x=\int_\Kx (\rho_h^n-\rho_i) w d\bx, $$
Note that \eqref{e0} is satisfied if the initial electric field is obtained exactly.

This also holds for DG methods with time integrators {\bf Scheme-2}$(\triangle t)$, {\bf Scheme-3}$(\triangle t)$, {\bf Scheme-2F}$(\triangle t)$, {\bf Scheme-3F}$(\triangle t)$.
\end{thm}

\emph{Proof.}  Take $g \in \mW_h^k \subset \mG_h^k$ in \eqref{dscheme1:3}, and sum over all $\Kv$,
\beq
\label{eqcontsc}
 \int_\Kx \frac{\rho_h^{n+1}-\rho_h^n}{\triangle t}  g d\bx
- \int_\Kx \bJ_h^{n+1/2}\cdot\nabla_\bx g d\bx
\quad+\int_{\df\Kx} \widehat{\bJ_h^{n+1/2} \cdot \bn_x} g ds_x  =0~,
\eeq
where $\rho_h^n(x, t)=\int_\Ov f_h^n d\bv \in  \mW_h^k$, and the flux is the central flux.
$$\widehat{\bJ_h^{n+1/2} \cdot \bn_x}:= \{\bJ_h^{n+1/2}\}_x \cdot\bn_x$$

We notice that \eqref{eqcontsc} is the DG scheme for the charge continuity equation $\rho_t+\nabla_\bx \cdot \bJ=0$ with central flux. Now assume that for any $w \in \mW_h^k$,
$$-\int_\Kx \bE_h^n \cdot \nabla_\bx w d\bx +\int_{\df\Kx} \widehat{\bE_h^n \cdot \bn_x} w ds_x=\int_\Kx (\rho_h^n-\rho_i) w d\bx, $$
where $$\widehat{\bE_h^n \cdot \bn_x}:= \{\bE_h^n\}_x \cdot\bn_x$$
is the central flux.
Then by \eqref{dscheme1:2} and \eqref{eqcontsc}
\begin{eqnarray*}
&&-\int_\Kx \bE_h^{n+1} \cdot \nabla_\bx w d\bx +\int_{\df\Kx} \widehat{\bE_h^{n+1} \cdot \bn_x} w ds_x \\
&&=-\int_\Kx \bE_h^n \cdot \nabla_\bx w d\bx +\int_{\df\Kx} \widehat{\bE_h^n \cdot \bn_x} w ds_x-\int_\Kx (\bE_h^{n+1}-\bE_h^n) \cdot \nabla_\bx w d\bx +\int_{\df\Kx} \widehat{(\bE_h^{n+1}-\bE_h^n) \cdot \bn_x} w ds_x \\
&&=\int_\Kx (\rho_h^n-\rho_i) w d\bx+\triangle t (\int_\Kx \bJ_h^{n+1/2} \cdot \nabla_\bx w d\bx -\int_{\df\Kx} \widehat{\bJ_h^{n+1/2}\cdot \bn_x} w ds_x )\\
&&=\int_\Kx (\rho_h^{n+1}-\rho_i) w d\bx
\end{eqnarray*}
By induction, we are done.  $\Box$

The theorem above shows that our scheme computes an electric field that is consistent with the Poisson equation, and this guarantees a physical relevant solution. Unfortunately, for the upwind flux, the derivation cannot go through, because we cannot find a simple DG scheme for the continuity equation as in \eqref{eqcontsc}.

We  also establish fully discrete $L^2$ stability for the implicit schemes. 
\begin{thm}[{$L^2$ stability}]
The  DG methods with time integrators {\bf Scheme-2}$(\triangle t)$, {\bf Scheme-3}$(\triangle t)$, {\bf Scheme-2F}$(\triangle t)$, {\bf Scheme-3F}$(\triangle t)$ satisfy
$$\int_\Ox \int_\Ov |f_h^\na|^2 d\bv d\bx =\int_\Ox \int_\Ov |f_h^n|^2  d\bv d\bx $$
for central flux, and
$$\int_\Ox \int_\Ov |f_h^\na|^2 d\bv d\bx \leq \int_\Ox \int_\Ov |f_h^n|^2  d\bv d\bx $$
for upwind flux.
\end{thm}
\emph{Proof.}  The proof is straightforward by taking the test function $g=\frac{1}{2}(f^n+f^{n+1})$ and is omitted.  $\Box$

\subsubsection{Splitted schemes and their properties}

In this subsection, we would like to design  fully discrete implicit schemes with the split-in-time integrators {\bf Scheme-4}$(\triangle t)$, {\bf Scheme-4F}$(\triangle t)$. The key idea is to solve each splitted equation in their respective reduced dimensions.
Let's introduce some notations for the description of the scheme. We look for
$f_h \in \mS_h^{k}$, for which
 we can pick a few nodal points to represent the degree of freedom for that element \cite{hesthaven2007nodal}. Suppose the nodes in $\Kx$ and $\Kv$ are  $\bx_\Kx^{(l)}$, $\bv_\Kv^{(m)}$, $l=1, \ldots, dof(k1)$, $m=1, \ldots, dof(k2)$,  respectively, then any $g \in \mS_h^{k}$ can be uniquely represented as $g=\sum_{l,m} g(\bx_\Kx^{(l)}, \bv_\Kv^{(m)}) L_x^{(l)}(\bx) L_v^{(m)}(\bv)$ on $K$, where $L_x^{(l)}(\bx), L_v^{(m)}(\bv)$ denote the $l$-th and $m$-th Lagrangian interpolating polynomials in $\Kx$ and $\Kv$, respectively.

Under this setting, the equations for $f$ in the splitted equations (a), (b)  can be solved in reduced dimensions. For example,  equation (a), we can fix a nodal point in $\bv$, say $\bv_\Kv^{(m)}$, then solve $\partial_t f (\bv_\Kv^{(m)})  +  \bv_\Kv^{(m)}\cdot \nabla_{\bx} f  (\bv_\Kv^{(m)})= 0$ by a DG methods in the $\bx$ direction. We can use the time integrator discussed in the previous subsection,  and get a update of point values at $f(\bx_\Kx^{(l)}, \bv_\Kv^{(m)})$ for all $\Kx, l$.

The idea is similar for equation (b). We can fix a nodal point in $\bx$, say $\bx_\Kx^{(l)}$, then solve
\[ \left\{
  \begin{array}{l}
\partial_t f (\bx_\Kx^{(l)})+   \bE(\bx_\Kx^{(l)})  \cdot \nabla_\bv f  (\bx_\Kx^{(l)}) = 0~, \\[3mm]
\partial_t \bE (\bx_\Kx^{(l)}) = -\bJ (\bx_\Kx^{(l)}),
  \end{array} \right.
  \]
 in $\bv$ direction, and get a update of point values at $f(\bx_\Kx^{(l)}, \bv_\Kv^{(m)})$ for all $\Kv, m$.


For simplicity of discussion, below we will describe in detail the scheme in a 1D1V setting.
For one-dimensional problems, we use a mesh that is a tensor product of grids in the $x$ and $v$ directions, and the domain $\Omega$ is  partitioned as follows:
\begin{equation*}
\label{2dcell1} 0=x_{\frac{1}{2}}<x_{\frac{3}{2}}< \ldots
<x_{N_x+\frac{1}{2}}=L , \qquad -V_c=v_{\frac{1}{2}}<v_{\frac{3}{2}}<
\ldots <v_{N_v+\frac{1}{2}}=V_c,
\end{equation*}
where $V_c$ is chosen appropriately large to guarantee $f(x, v, t)=0$ for $|v| \geq V_c$.   The grid is defined as
\begin{eqnarray*}
\label{2dcell2} && K_{i,j}=[x_{i-\frac{1}{2}},x_{i+\frac{1}{2}}]
\times [v_{j-\frac{1}{2}},v_{j+\frac{1}{2}}] , \nonumber \\
&&K_{x,i}=[x_{i-1/2},x_{i+1/2}],
\quad K_{v,j}=[v_{j-1/2},v_{j+1/2}]\,, \     \quad i=1,\ldots N_x, \quad j=1,\ldots N_v ,
\end{eqnarray*}
Let  $\triangle x_i= x_{i+1/2}-x_{i-1/2}$, $\triangle v_j=v_{j+1/2}-v_{j-1/2}$ be the length of each interval. $x_i^{(l)}, l=1, \ldots, k+1$ be the $(k+1)$ Gauss quadrature points on $K_{x,i}$ and $v_j^{(m)}, m=1, \ldots, k+1$ be the $(k+1)$ Gauss quadrature points on $K_{v,j}$. Now we are ready to describe our scheme.
\bigskip

\underline{Algorithm {\bf Scheme-a}$(\triangle t)$}
\medskip

To solve from $t^n$ to $t^\na$
 \[
  \textrm{ (a)} \left\{
  \begin{array}{l}
\partial_t f  +  v f_x   = 0~, \\
\partial_t E = 0,
  \end{array} \right.
\]
\begin{enumerate}
\item For each $j=1,\ldots N_v , m=1, \ldots, k+1$, we seek $g_j^{(m)}(x) \in \mW_h^k$, such that
\begin{eqnarray}
 &&   \int_{J_i} \frac{g_j^{(m)}(x)-f_h^n(x, v_j^{(m)} )}{\triangle t} \varphi_h \, dx - \int_{J_i} v_j^{(m)}\frac{g_j^{(m)}(x)+f_h^n(x, v_j^{(m)} )}{2} (\varphi_h)_x \, dx \label{schemeas} \\
 &&+ \widehat{v_j^{(m)} \frac{g_j^{(m)}(x_{i+\frac{1}{2}})+f_h^n(x_{i+\frac{1}{2}}, v_j^{(m)} )}{2} } (\varphi_h)_{i+\frac{1}{2}}^- \,
    -\widehat{v_j^{(m)} \frac{g_j^{(m)}(x_{i-\frac{1}{2}})+f_h^n(x_{i-\frac{1}{2}}, v_j^{(m)} )}{2} } (\varphi_h)_{i-\frac{1}{2}}^+ =0 \notag
\end{eqnarray}
holds for any test function $\varphi_h(x,t) \in \mW_h^k$.
\item Let $f_h^{n+1}$ be the unique polynomial in $\mS_h^k$, such that $f_h^\na(x_i^{(l)}, v_j^{(m)})=g_j^{(m)}(x_i^{(l)}), \, \forall i, j, l, m$.
\end{enumerate}

\bigskip

\underline{Algorithm {\bf Scheme-b}$(\triangle t)$}
\medskip

To solve from $t^n$ to $t^\na$
\[
  \textrm{(b)} \left\{
  \begin{array}{l}
\partial_t f  +   Ef_v  = 0~, \\
\partial_t E = -J,
  \end{array} \right.
\]

\begin{enumerate}
\item For each $i=1,\ldots N_x , l=1, \ldots, k+1$, we seek $g_i^{(l)}(v) \in \mZ_h^k$ and $E_i^{(l)}$, such that
\begin{eqnarray}
  &&\int_{K_j} \frac{g_i^{(l)}(v)-f_h^n(x_i^{(l)}, v)}{\triangle t} \varphi_h \,  dv - \int_{K_j} \frac{E_h^n(x_i^{(l)})+E_i^{(l)}}{2}  \frac{g_i^{(l)}(v)+f_h^n(x_i^{(l)}, v)}{2}  (\varphi_h)_v \, dv \notag\\
   && + \frac{E_h^n(x_i^{(l)})+E_i^{(l)}}{2} \widehat{\frac{g_i^{(l)}(v_{j+\frac{1}{2}})+f_h^n(x_i^{(l)}, v_{j+\frac{1}{2}})}{2}} (\varphi_h)^-_{ j+\frac{1}{2}} \notag\\
    &&-\frac{E_h^n(x_i^{(l)})+E_i^{(l)}}{2} \widehat{\frac{g_i^{(l)}(v_{j-\frac{1}{2}})+f_h^n(x_i^{(l)}, v_{j-\frac{1}{2}})}{2}} (\varphi_h)^+_{ j-\frac{1}{2}} =0 \label{scheme3ss}\\
   &&\frac{E_i^{(l)}-E_h^n(x_i^{(l)})}{\triangle t}=  -\frac{1}{2}(J_h^n(x_i^{(l)})+J_i^{(l)}), \; \textrm{where} \, J_h^n(x)=\int_\Ov f_h^n(x, v) v dv, J_i^{(l)}=\int_\Ov g_i^{(l)}(v) v dv \notag
\end{eqnarray}
holds for any test function $\varphi_h(x,t) \in \mZ_h^k$.
\item Let $f_h^{n+1}$ be the unique polynomial in $\mS_h^k$, such that $f_h^\na(x_i^{(l)}, v_j^{(m)})=g_i^{(l)}(x_j^{(m)}),\, \forall i, j, l, m$.
Let $E_h^{n+1}$ be the unique polynomial in $\mW_h^k$, such that $E_h^\na(x_i^{(l)})=E_i^{(l)}, \, \forall i, l$.
\end{enumerate}

The flux terms in the algorithms above could be either upwind or central flux. Finally, we recall
$\textnormal{\bf Scheme-4}(\triangle t)=\textnormal{\bf Scheme-a}(\triangle t/2) \textnormal{\bf Scheme-b}(\triangle t) \textnormal{\bf Scheme-a}(\triangle t/2)$
 and $\textnormal{\bf Scheme-4F}(\triangle t)=\textnormal{\bf Scheme-4}(\beta_1 \triangle t) \textnormal{\bf Scheme-4}(\beta_2 \triangle t) \textnormal{\bf Scheme-4}(\beta_3 \triangle t)$.

\bigskip

Next, we will discuss the conservation properties of the methods above.

\begin{thm}[{Total particle number conservation}]
The DG schemes with  {\bf Scheme-4}$(\triangle t)$, {\bf Scheme-4F}$(\triangle t)$
  preserve the total particle number of the system, i.e.
$$\int_\Ox \int_\Ov f_h^\na dv dx =\int_\Ox \int_\Ov f_h^n  dv dx .$$
\end{thm}

\emph{Proof.}  We only need to show that for each of the operators $\textnormal{\bf Scheme-a}(\triangle t)$ and $\textnormal{\bf Scheme-b}(\triangle t)$. For $\textnormal{\bf Scheme-a}(\triangle t)$,
let $\varphi_h=1$ in \eqref{schemeas}, and sum over all element $K_{i,x}$, we get
$$\int_\Ox   g_j^{(m)}(x) dx =\int_\Ox f_h^n(x, v_j^{(m)}) dx .$$
Therefore for any $j, m$,
 $$\int_\Ox  f_h^\na(x, v_j^{(m)})  dx =\int_\Ox f_h^n(x, v_j^{(m)}) dx ,$$
 and because the (k+1)-point Gauss quadrature formula is exact for polynomial with degree less than $2k+2$,
\begin{eqnarray*}
\int_\Ox \int_\Ov f_h^\na dv dx =\sum_j \sum_m w_m \int_\Ox  f_h^\na(x, v_j^{(m)})  dx \,\triangle v_j\\
=\sum_j \sum_m w_m \int_\Ox  f_h^n(x, v_j^{(m)})  dx \,\triangle v_j=\int_\Ox \int_\Ov f_h^n  dv dx ,
\end{eqnarray*}
 where $w_m$ is the corresponding Gauss quadrature weights. The proof is similar for $\textnormal{\bf Scheme-b}(\triangle t)$ and is omitted.
 $\Box$

\begin{thm}[{Total energy conservation}]
If $k \geq 2$, the   DG schemes with  $\textnormal{\bf Scheme-4}(\triangle t)$, $\textnormal{\bf Scheme-4F}(\triangle t)$ preserve the discrete total energy $TE_n=TE_{n+1}$, where
$$2 \,(TE_n)=\int_\Ox \int_\Ov f_h^n |v|^2 dv dx +\int_\Ox |E_h^n|^2 dx.$$
\end{thm}

\emph{Proof.}  We  need to show that each of the operators $\textnormal{\bf Scheme-a}(\triangle t)$ and $\textnormal{\bf Scheme-b}(\triangle t)$ preserves the total energy. For $\textnormal{\bf Scheme-a}(\triangle t)$,
similar to the previous proof, we obtain that for any $j, m$,
 $$\int_\Ox  f_h^\na(x, v_j^{(m)})  dx =\int_\Ox f_h^n(x, v_j^{(m)}) dx ,$$
 and because the (k+1) Gauss quadrature formula is exact for polynomial with degree less than $2k+2$,
\begin{eqnarray*}
\int_\Ox \int_\Ov f_h^\na |v|^2dv dx =\sum_j \sum_m w_m \int_\Ox  f_h^\na(x, v_j^{(m)})  dx \, |v_j^{(m)}|^2 \,\triangle v_j\\
=\sum_j \sum_m w_m \int_\Ox  f_h^n(x, v_j^{(m)})  dx \, |v_j^{(m)}|^2 \,\triangle v_j=\int_\Ox \int_\Ov f_h^n   |v|^2dv dx ,
\end{eqnarray*}
because $\int_\Ox f_h(x, v) dx |v|^2$ is a polynomial in $v$ that is at most degree $k+2$, since $k \geq 2$, we know $k+2 \leq 2k+1$, and the Gauss quadrature formula is exact.

As for $\textnormal{\bf Scheme-b}(\triangle t)$,  because $|E_h|^2$ is a polynomial of degree at most $2k$, therefore the Gauss quadrature is exact, and using the scheme
\begin{eqnarray*}
\int_\Ox |E_h^\na|^2-|E_h^n|^2 dx&=&\sum_i \sum_l w_l \, \left ( (E_h^\na(x_i^{(l)}))^2- (E_h^n(x_i^{(l)}))^2 \right) \triangle x_i\\
&=&\sum_i \sum_l w_l \, \left ( (E_i^{(l)}))^2- (E_h^n(x_i^{(l)}))^2 \right) \triangle x_i\\
&=&\sum_i \sum_l w_l \, \left ( E_i^{(l)})+ E_h^n(x_i^{(l)}) \right)  \left (  -\frac{\triangle t}{2}(J_h^n(x_i^{(l)})+J_i^{(l)}) \right)  \triangle x_i
\end{eqnarray*}
Because $k \geq 2$, we can take $\varphi_h=v^2$ in \eqref{scheme3ss}, and
\begin{eqnarray*}
&&\int_\Ox \int_\Ov f_h^\na |v|^2 dv dx-\int_\Ox \int_\Ov f_h^n |v|^2 dv dx\\
&=&\sum_i \sum_l w_l \, \left ( \int_\Ov f_h^\na (x_i^{(l)}, v)  v^2 dv-\int_\Ov f_h^n (x_i^{(l)}, v)  v^2 dv \right) \triangle x_i\\
&=&\sum_i \sum_l w_l \, \left ( \int_\Ov (g_i^{(l)}-f_h^n (x_i^{(l)}, v))  v^2 dv \right) \triangle x_i\\
&=&\sum_i \sum_l w_l \, \left ( \triangle t  \int_{\Ov} \frac{E_h^n(x_i^{(l)})+E_i^{(l)}}{2}  (g_i^{(l)}(v)+f_h^n(x_i^{(l)}, v))v \, dv \right) \triangle x_i\\
&=&\sum_i \sum_l w_l \, \left ( \triangle t   \frac{E_h^n(x_i^{(l)})+E_i^{(l)}}{2}  (J_h^n(x_i^{(l)})+J_i^{(l)})\right) \triangle x_i
\end{eqnarray*}
Therefore, putting the results for $\textnormal{\bf Scheme-a}(\triangle t)$ and $\textnormal{\bf Scheme-b}(\triangle t)$ together, we have proved
$$\int_\Ox \int_\Ov f_h^\na |v|^2 dv dx +\int_\Ox |E_h^\na|^2 dx=\int_\Ox \int_\Ov f_h^n |v|^2 dv dx +\int_\Ox |E_h^n|^2 dx.$$
 $\Box$

\medskip

\begin{thm}[{$L^2$ stability}]
The DG schemes with  $\textnormal{\bf Scheme-4}(\triangle t)$, $\textnormal{\bf Scheme-4F}(\triangle t)$
satisfy
$$\int_\Ox \int_\Ov |f_h^\na|^2 dv dx =\int_\Ox \int_\Ov |f_h^n|^2  dv dx $$
for central flux, and
$$\int_\Ox \int_\Ov |f_h^\na|^2 dv dx \leq \int_\Ox \int_\Ov |f_h^n|^2  dv dx $$
for upwind flux.
\end{thm}
\emph{Proof.}   We only need to show that each of the operators $\textnormal{\bf Scheme-a}(\triangle t)$ and $\textnormal{\bf Scheme-b}(\triangle t)$. For  $\textnormal{\bf Scheme-a}(\triangle t)$,
let $\varphi_h=g_j^{(m)}(x)+ f_h^n(x, v_j^{(m)}) $ in \eqref{schemeas}, and sum over all element $K_{x,i}$, we get
$$\int_\Ox  | g_j^{(m)}(x)|^2 dx =\int_\Ox |f_h^n(x, v_j^{(m)})|^2 dx $$
for central flux and
$$\int_\Ox  | g_j^{(m)}(x)|^2 dx \leq \int_\Ox |f_h^n(x, v_j^{(m)})|^2 dx $$
for upwind flux. Therefore for any $j, m$,
 $$\int_\Ox  |f_h^\na(x, v_j^{(m)})|^2  dx =\int_\Ox |f_h^n(x, v_j^{(m)})|^2 dx ,$$
 for central flux and
  $$\int_\Ox  |f_h^\na(x, v_j^{(m)})|^2  dx \leq \int_\Ox |f_h^n(x, v_j^{(m)})|^2 dx ,$$
  for upwind flux.
 and because the (k+1) Gauss quadrature formula is exact for polynomial with degree less than $2k+2$,
\begin{eqnarray*}
\int_\Ox \int_\Ov |f_h^\na|^2 dv dx =\sum_j \sum_m w_m \int_\Ox  |f_h^\na(x, v_j^{(m)})|^2  dx \,\triangle v_j\\
= (\textrm{or} \leq) \sum_j \sum_m w_m \int_\Ox  |f_h^n(x, v_j^{(m)})|^2  dx \,\triangle v_j=\int_\Ox \int_\Ov |f_h^n|^2  dv dx ,
\end{eqnarray*}
 where $w_m$ is the corresponding Gauss quadrature weights. The proof is similar for $\textnormal{\bf Scheme-b}(\triangle t)$ and is omitted.
 $\Box$


In summary, this  type of operator splitting enable us to solve the VA system by solving two essentially decoupled lower-dimensional equations, and still maintain energy conservation.  A fully implicit method is possible, and it does not require to invert a coupled nonlinear systems in the whole $(\bx, \bv)$ space; and very large time step is allowed in this framework when implicit schemes are used.
In numerical computation, a
Jacobian-free Newton-Krylov solver \cite{knoll2004jacobian} is used to compute the nonlinear systems resulting from implicit scheme \eqref{scheme3ss}.


\section{Numerical Results}
\label{sec:numerical}

In this section, we tests   the proposed methods  by the following benchmark examples:
\begin{exa}
(Landau damping)  
In this case, the initial condition is given by
  \begin{align}
    \label{landau}
      f_0(x,v)=f_M(v)(1+A \cos(\kappa x)), \quad x\in[0, L],\; v\in[-V_c, V_c],
   \end{align}
   with $A = 0.5, \; \kappa = 0.5, L = 4\pi,\; V_c = 8$, and
    $f_M = \frac{1}{\sqrt{2\pi}}e^{-v^2/2}.$
\end{exa}

\begin{exa}
(Two-stream instability)
In this case, the initial condition is given by
\begin{align}
    \label{twostream}
      f_0(x,v)=f_{TS}(v)(1+A \cos(\kappa x)), \quad x\in[0, L],\; v\in[-V_c, V_c],
   \end{align}
   with $A = 0.05, \; \kappa = 0.5, L = 4\pi,\; V_c = 8$, and
    $f_{TS} = \frac{1}{\sqrt{2\pi}}v^2e^{-v^2/2}.$
\end{exa}

\begin{exa}
(Bump-on-tail instability \cite{banks2010new,hittinger2013block})
In this case, the initial condition is given by
\begin{align}
    \label{bump}
      f_0(x,v)=f_{BT}(v)(1+A \cos(\kappa x)), \quad x\in[0, L],\; v\in[-V_c, V_c],
   \end{align}
   with $A = 0.04, \; \kappa = 0.3, L = 20\pi/3,\; V_c = 13$; and
    $$f_{BT} = n_p exp(-v^2/2)+n_b exp\left(-\frac{|v-u|^2}{2v_t^2}\right)$$
    whose parameters are 
      $$n_p = \frac{9}{10\sqrt{2\pi}}, \quad n_b = \frac{2}{10\sqrt{2\pi}},\quad u= 4.5, \quad v_t = 0.5.$$
\end{exa}

Note that in the examples above, we have taken $V_c$  to be larger than the common values in the literature to eliminate the boundary effects and to accurately reflect the conservation properties of our methods. For $\textnormal{\bf Scheme-4}$, $\textnormal{\bf Scheme-4F}$, we use KINSOL from SUNDIALS \cite{hindmarsh2005sundials} to solve the nonlinear algebraic systems resulted from the discretization of equation (b), and we  set the tolerance parameter to be $\epsilon_{tol}=10^{-12}$. In all the runs below, for simplicity, we use uniform meshes in $x$ and $v$ directions, while we note that nonuniform mesh can also be used under the DG framework.

\subsection{Accuracy tests}
In this subsection, we use two-stream instability  to test the orders of accuracy of the proposed schemes. Because the VA system is
time reversible, if we let
$f(x,v,0),\;E(x,0)$ be the the initial conditions of the VA system,   $f(x,v,T),\,E(x,T)$ be the solution at $t=T$, and
  we enforce $f(x,-v,T),\,E(x,T)$ be the initial condition for the VA system at $t=0$, then at $t=T$, we would recover
$f(x,-v,0),\;E(x,0)$. This provides a way to measure the errors of our schemes. In Tables \ref{table_explicit2nd} to \ref{table_strang4th},  we run the VA system to $T=5$ and then back to $T=10$, and compare it with the initial conditions.

%
%
%
 
 Table \ref{table_explicit2nd} lists the $L^2$ error and
order for the explicit scheme $\textnormal{\bf Scheme-1}$. To match the accuracy of the temporal and spatial discretizations, we take   $\Delta t\sim \min(\Delta x, \Delta v)$ for space $\mG_h^{1}$, and $\Delta t\sim \min(\Delta x, \Delta v)^{3/2}$ for $\mG_h^{2}$, and $\Delta t\sim \min(\Delta x, \Delta v)^2$ for $\mG_h^{3}$. This is also done for Table \ref{table_splittime2nd} and  \ref{table_strang2nd}. Because of the stability restriction of the explicit scheme, we take $CFL = 0.3$ for $\mG_h^{1}$,  and the coefficient to be $0.13, 0.05$ for $\mG_h^{2}$ and $\mG_h^{3}$, respectively. From the table, we can see that for all three polynomial spaces, we obtain the optimal $(k+1)$-th order for $f$, while the convergence order of $E$ is higher.

 Table \ref{table_splittime2nd}-\ref{table_strang4th} are for the implicit methods, and we take all time coefficient to be 2.
Table \ref{table_splittime2nd} and \ref{table_splittime4th} list the $L^2$ error and order of schemes $\textnormal{\bf Scheme-4}$ and $\textnormal{\bf Scheme-4F}$. The errors behave similarly as the explicit scheme $\textnormal{\bf Scheme-1}$, except for the $E$ component of the $\textnormal{\bf Scheme-4}$ scheme with $160 \times 160$ mesh.  The order for this mesh is decreased, because the tolerance parameter $\epsilon_{tol}$ for the Newton-Krylov solve has polluted the error in this calculation. 
Table \ref{table_strang2nd} and \ref{table_strang4th}  list the $L^2$ error and order of the schemes $\textnormal{\bf Scheme-3F}$ and $\textnormal{\bf Scheme-3F}$.  We also achieve the optimal $(k+1)$-th order for $f$ and higher order convergence for $E$.

 \begin{table}[!htbp]
\caption{Time discretization $\textnormal{\bf Scheme-1}$, with DG methods using the indicated polynomial space and the upwind flux. $L^2$ error and order.}
\label{table_explicit2nd}
\smallskip
\begin{tabular}{|c|c|c|c|c|c|c|c|c|}
\hline
\multicolumn{2}{|c|}{\multirow{2}*{\backslashbox{Space}{Mesh }}}&$20\times 20$&\multicolumn{2}{|c|}{$40\times 40$}&\multicolumn{2}{c}{$80\times 80$}&\multicolumn{2}{|c|}{$160\times 160$}\\
\cline{3-9}
 \multicolumn{2}{|c|}{} & Error& Error & Order& Error & Order&Error & Order\\
\hline
\multirow{2}*{$\mG_h^{1}$}
&f&	7.32E-03&	1.84E-03&	1.99	&5.05E-04&	1.87	&1.40E-04&	1.85\\
\cline{2-9}
&E&	7.51E-03&	1.22E-03&	2.62&	1.60E-04&	2.93&	1.95E-05&	3.03\\
\hline
\multirow{2}*{$\mG_h^{2}$}
&f&	1.21E-03&	2.09E-04&	2.53&	3.16E-05	&2.73	&3.26E-06	&3.28\\
\cline{2-9}
&E&	3.42E-04	&1.56E-05&	4.46	&7.97E-07&	4.29	&7.89E-08&	3.34\\
\hline
\multirow{2}*{$\mG_h^{3}$}
&f&	2.43E-04&	1.85E-05&	3.72&	1.40E-06&	3.72&	7.68E-08	&4.19\\
\cline{2-9}
&E&	2.17E-05&	4.09E-07&	5.73&	5.44E-09&	6.23	&7.78E-11&	6.13\\
\hline
\end{tabular}
\end{table}


\begin{table}[!htbp]
\caption{Time discretization $\textnormal{\bf Scheme-4}$, with DG methods using the indicated polynomial space and the upwind flux. $L^2$ error and order.}
\label{table_splittime2nd}
\smallskip
\begin{tabular}{|c|c|c|c|c|c|c|c|c|}
\hline
\multicolumn{2}{|c|}{\multirow{2}*{\backslashbox{Space}{Mesh }}}&$20\times 20$&\multicolumn{2}{|c|}{$40\times 40$}&\multicolumn{2}{c}{$80\times 80$}&\multicolumn{2}{|c|}{$160\times 160$}\\
\cline{3-9}
 \multicolumn{2}{|c|}{} & Error& Error & Order& Error & Order&Error & Order\\
\hline
\multirow{2}*{$\mS_h^{1}, \mW_h^1$}&
f&4.03E-03&	1.08E-03&	1.90&	3.46E-04	&1.64	&1.05E-04	&1.72\\
\cline{2-9}
&E	&1.16E-03&	1.68E-04	&2.79&	3.13E-05&	2.42&	5.87E-06	&2.41\\
\hline
\multirow{2}*{$\mS_h^2,\mW_h^2$}
&f&	8.03E-04	&1.65E-04&	2.28&	2.82E-05&	2.55&	2.73E-06&	3.37\\
\cline{2-9}
&E&	6.11E-05	&3.49E-06&	4.13&	2.31E-07&	3.92&	1.39E-08&	4.05\\
\hline
\multirow{2}*{$\mS_h^{3},\mW_h^3$}
&f&	1.54E-04&	1.22E-05&	3.66&	1.13E-06&	3.43&	6.64E-08&	4.09\\
\cline{2-9}
&E&	3.96E-06&	6.94E-08&	5.83&	9.26E-10&	6.23&	8.72E-10&	0.09\\\hline
\end{tabular}
\end{table}

\begin{table}[!htbp]
\caption{Time discretization $\textnormal{\bf Scheme-4F}$, with DG methods using the indicated polynomial space and the upwind flux. $L^2$ error and order. }
\label{table_splittime4th}
\smallskip
\begin{tabular}{|c|c|c|c|c|c|c|c|c|}
\hline
\multicolumn{2}{|c|}{\multirow{2}*{\backslashbox{Space}{Mesh }}}&$20\times 20$&\multicolumn{2}{|c|}{$40\times 40$}&\multicolumn{2}{c}{$80\times 80$}&\multicolumn{2}{|c|}{$160\times 160$}\\
\cline{3-9}
 \multicolumn{2}{|c|}{} & Error& Error & Order& Error & Order&Error & Order\\
\hline
\multirow{2}*{$\mS_h^{3},\mW_h^3$}
&f&	1.25E-04&	7.54E-06&	4.05&	5.12E-07&	3.88&	2.40E-08&4.42\\
\cline{2-9}
&E&	3.03E-06&	6.65E-08&	5.51&	4.15E-09&	4.00&	3.02E-10&3.78\\
\hline
\end{tabular}
\end{table}
\begin{table}[!htbp]
\centering
\caption{Time discretization $\textnormal{\bf Scheme-3}$, with DG methods using the indicated polynomial space and the upwind flux. $L^2$ error and order. }
\label{table_strang2nd}
\smallskip
\begin{tabular}{|c|c|c|c|c|c|c|}
\hline
\multicolumn{2}{|c|}{\multirow{2}*{\backslashbox{Space}{Mesh }}}&$20\times 20$&\multicolumn{2}{|c|}{$40\times 40$}&\multicolumn{2}{|c|}{$80\times 80$}\\
\cline{3-7}
 \multicolumn{2}{|c|}{} & Error& Error & Order& Error & Order\\
\hline
\multirow{2}*{$\mG_h^{1}$}
&f&	7.34E-03	&1.84E-03&	2.00	&5.05E-04&	1.87	\\
\cline{2-7}
&E&	7.33E-03	&1.21E-03&	2.60&	1.60E-04&	2.92\\
\hline
\multirow{2}*{$\mG_h^{2}$}
&f&	1.21E-03&	2.09E-04&	2.53&	3.16E-05&	2.73\\
\cline{2-7}
&E&	3.42E-04&	1.56E-05&	4.45&	7.98E-07&	4.29	\\
\hline
\multirow{2}*{$\mG_h^{3}$}
&f&	2.43E-04&	1.85E-05&	3.72&	1.40E-06&	3.72\\
\cline{2-7}
&E&	2.06E-05&	3.99E-07&	5.69&	8.21E-09	&5.60\\
\hline
\end{tabular}
\end{table}
	
\begin{table}[!htbp]
\centering
\caption{Time discretization $\textnormal{\bf Scheme-3F}$, with DG methods using the indicated polynomial space and the upwind flux. $L^2$ error and order. }
\label{table_strang4th}
\smallskip
\begin{tabular}{|c|c|c|c|c|c|c|}
\hline
\multicolumn{2}{|c|}{\multirow{2}*{\backslashbox{Space}{Mesh }}}&$20\times 20$&\multicolumn{2}{|c|}{$40\times 40$}&\multicolumn{2}{|c|}{$80\times 80$}\\
\cline{3-7}
 \multicolumn{2}{|c|}{} & Error& Error & Order& Error & Order\\
\hline
\multirow{2}*{$\mG_h^{3}$}
&f&	2.78E-04&	3.34E-05&	3.06&	2.25E-06&	3.89\\
\cline{2-7}
&E&	3.93E-05&	5.53E-07&	6.15&	2.35E-08&	4.56\\
\hline
\end{tabular}
\end{table}

\subsection{Conservation properties}
In this subsection, we will verify the conservation properties of the proposed methods $\textnormal{\bf Scheme-1}$, $\textnormal{\bf Scheme-4}$, and $\textnormal{\bf Scheme-4F}$. In particular, we test $\textnormal{\bf Scheme-1}$ with space $\mG_h^2$ and let  $CFL =0.13$, and denote it by $\textnormal{\bf Scheme-1}\textnormal{ and } P^2$ in the figures. We take $\textnormal{\bf Scheme-4}$ with space $\mS_h^2$ and we denote it by $\textnormal{\bf Scheme-4}\textnormal{ and } Q^2$, $\textnormal{\bf Scheme-4F}$ with space $\mS_h^3$ and denote it by $\textnormal{\bf Scheme-4} \textnormal { and } Q^2$. For those implicit runs on $100 \times 200$ mesh, we use $CFL=5$, except for Landau damping with  $\textnormal{\bf Scheme-4F} \textnormal{ and } Q^3$, for which we have taken $CFL=10$. To save computational time, we take $CFL=10$ and $\epsilon_{tol}=10^{-10}$ for figures computed by $\textnormal{\bf Scheme-3}$ 
and $\textnormal{\bf Scheme-3F}$. 

\begin{figure}[!htbp]

 \subfigure[Landau damping. Total particle number.]{\includegraphics[width=0.45\textwidth]{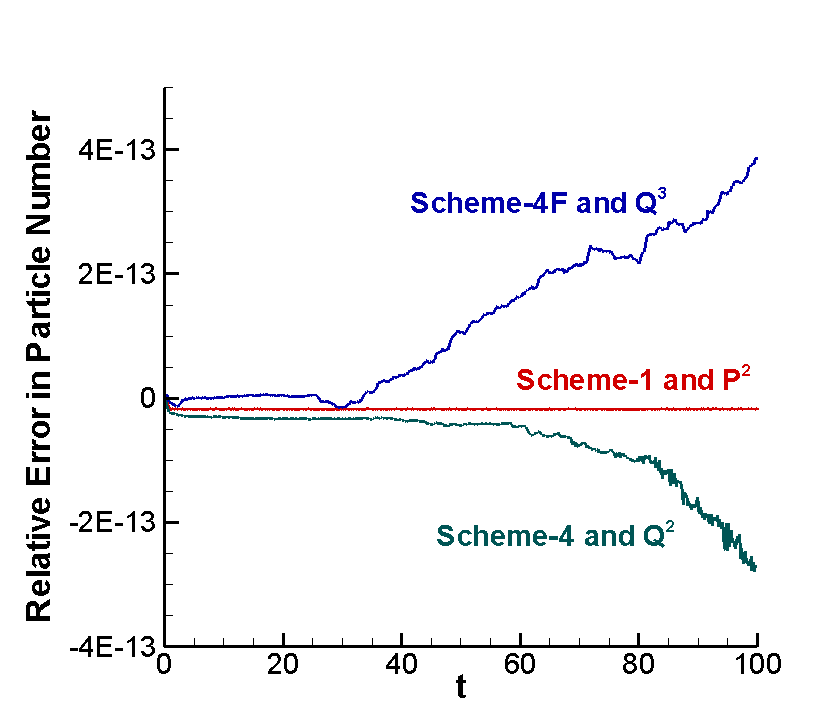}}
 \subfigure[Landau damping. Total energy.]{\includegraphics[width=0.45\textwidth]{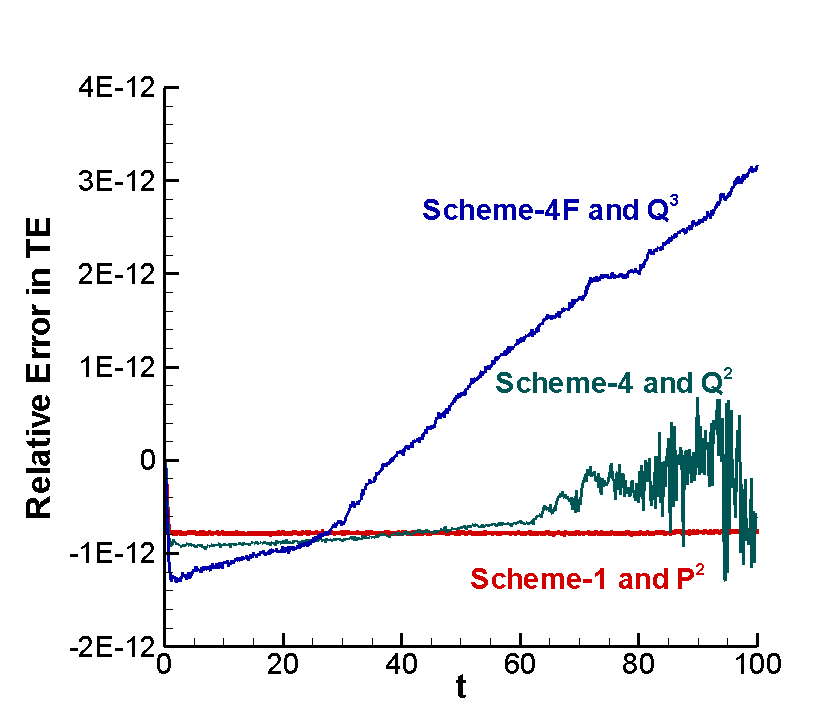}}\\
  \subfigure[Two-stream instability. Total particle number.]{\includegraphics[width=0.45\textwidth]{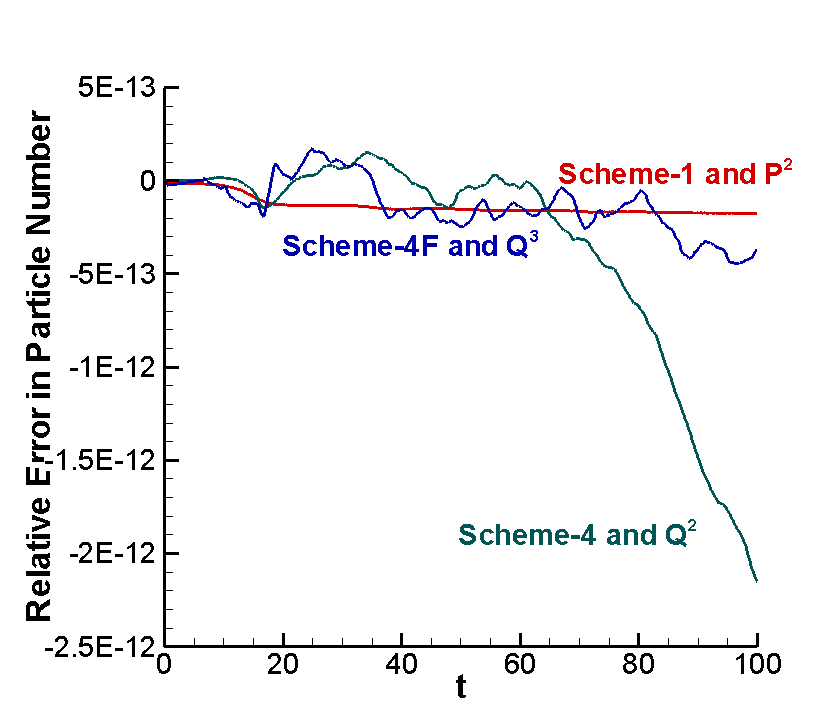}}
 \subfigure[Two-stream instability. Total energy.]{\includegraphics[width=0.45\textwidth]{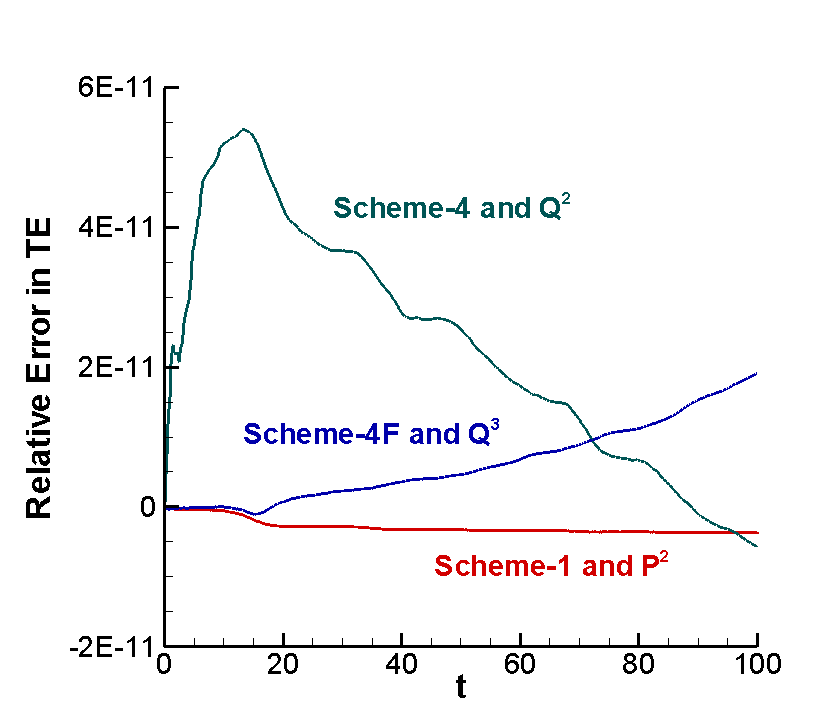}}\\
  \subfigure[Bump-on-tail instability. Total particle number.]{\includegraphics[width=0.45\textwidth]{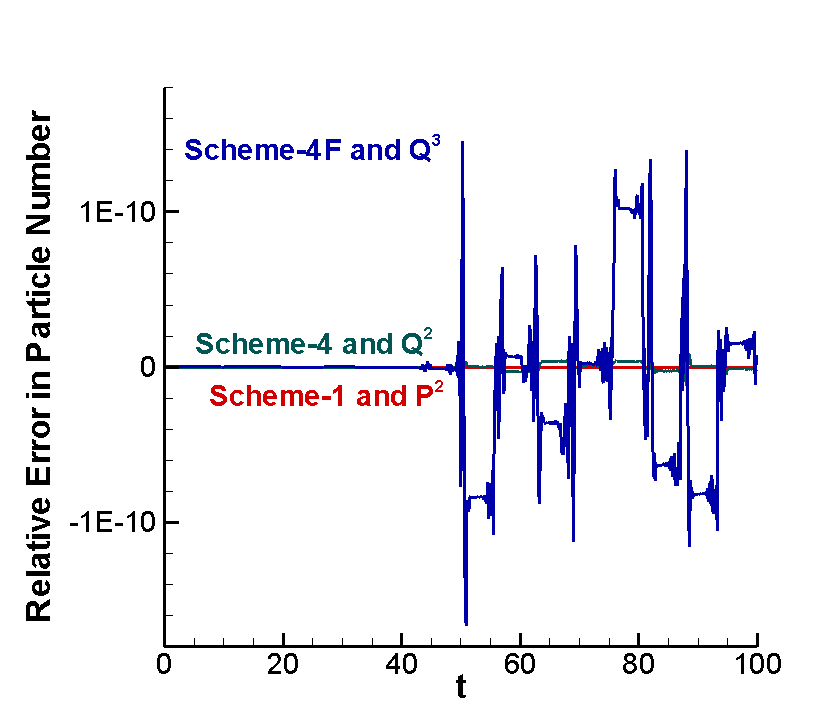}}
 \subfigure[Bump-on-tail instability. Total energy.]{\includegraphics[width=0.45\textwidth]{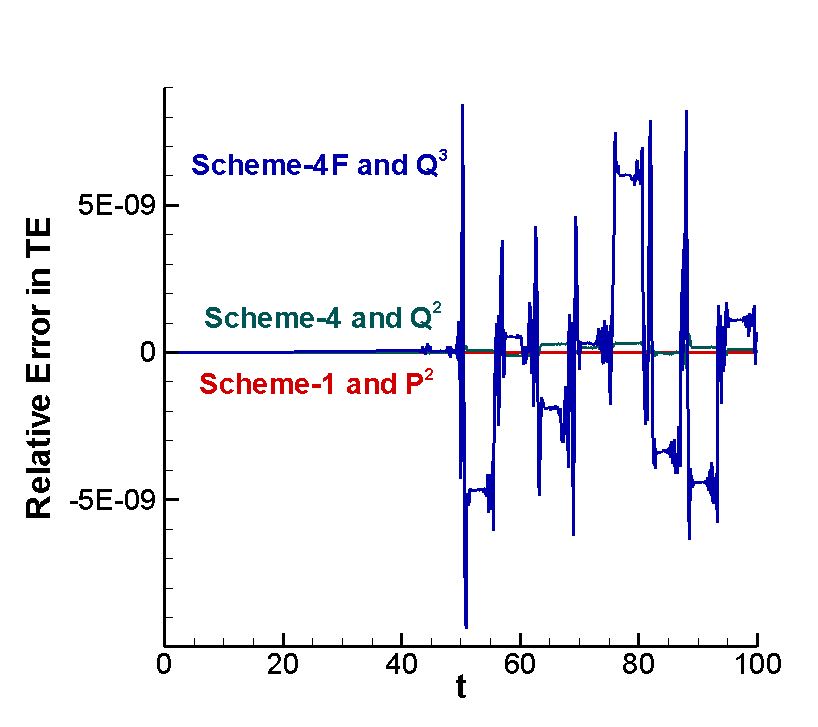}}
\caption{Evolution of the relative error in total particle number and total energy by the indicated methods. $100\times200$ mesh. Upwind flux.}
\label{figure_energy1}
\end{figure}

\begin{figure}[!htbp]
\subfigure[Landau damping.]{\includegraphics[width=0.45\textwidth]{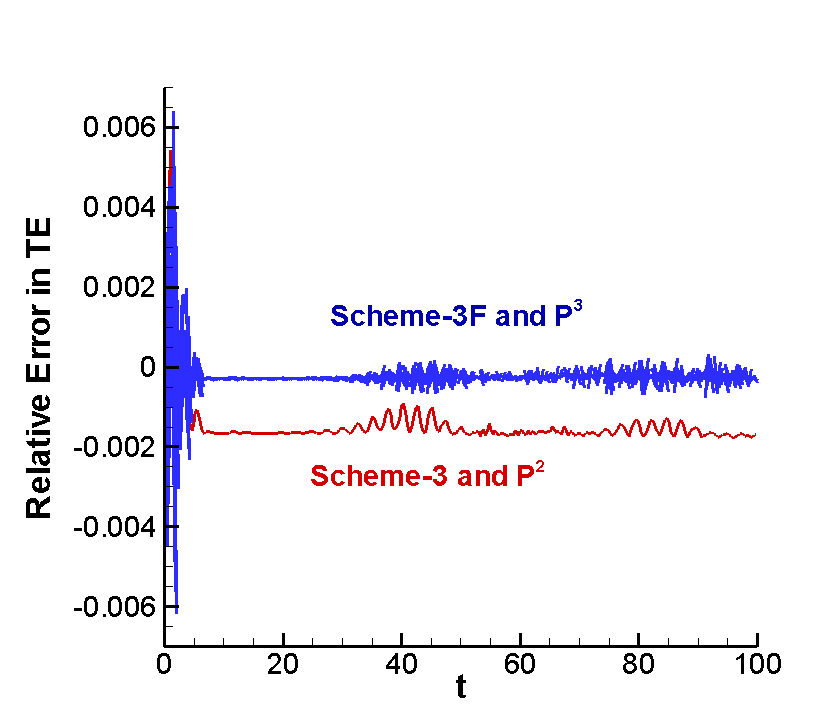}}
\subfigure[Two-stream instability.]{\includegraphics[width=0.45\textwidth]{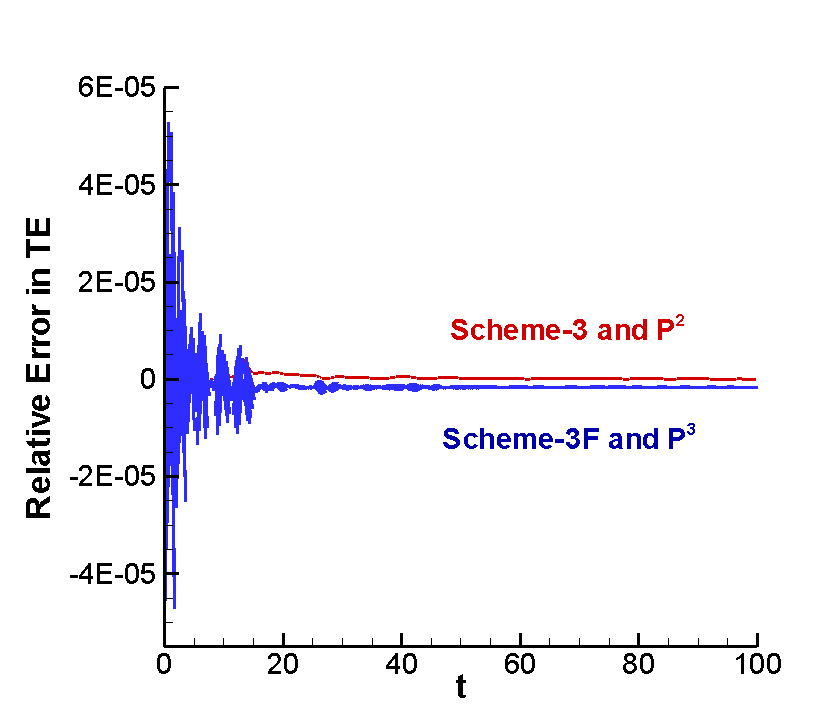}}
\\
\subfigure[Bump-on-tail instability]{\includegraphics[width=0.45\textwidth]{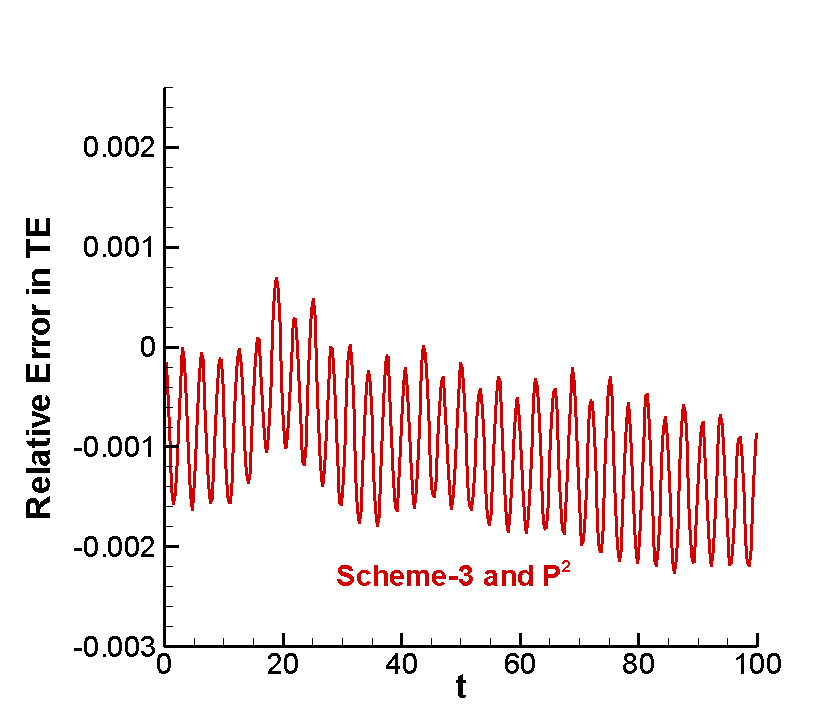}}
\subfigure[Bump-on-tail instability]{\includegraphics[width=0.45\textwidth]{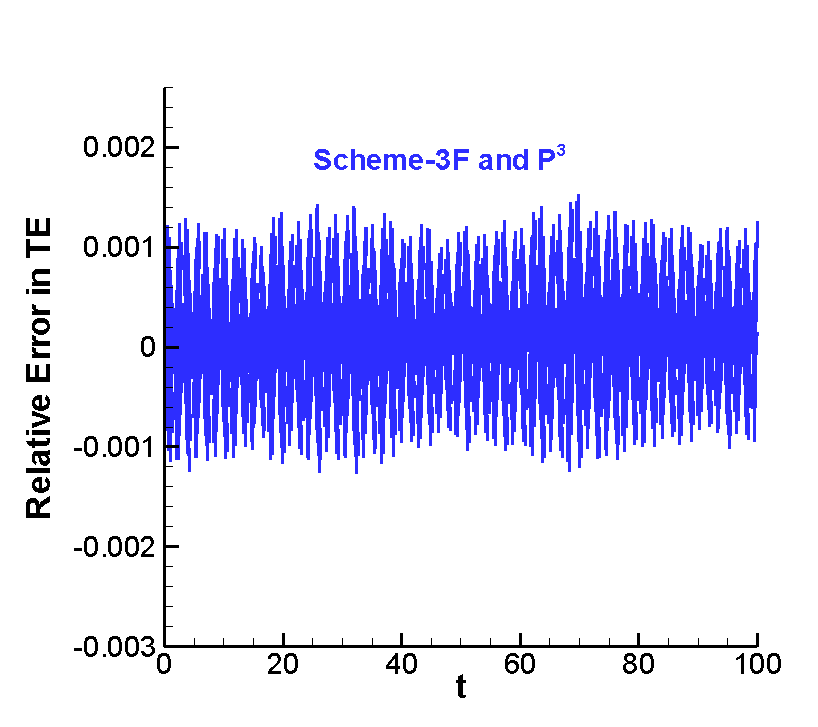}}
\caption{Evolution of the relative error in total energy by the indicated methods. $50\times 100$ mesh. Upwind flux. 
$ \epsilon_{tol}=1e-10$.
\label{figure_strang}}
\end{figure}

In Figure \ref{figure_energy1}, we plot the relative error of the total particle number and total energy for the three methods with the upwind flux. We  observe that the relative errors   stay small, below $10^{-11}$ for Landau damping, $10^{-10}$ for two-stream instability, and $10^{-8}$ for bump-on-tail instability. The conservation is especially good with the explicit method $\textnormal{\bf Scheme-1}$. As for the implicit schemes, $\textnormal{\bf Scheme-4}$ and $\textnormal{\bf Scheme-4F}$, the errors are larger mainly due to the error caused by   the Newton-Krylov solver, and is related to $\epsilon_{tol}$.
In Figure \ref{figure_strang}, we plot the relative error of the total energy obtained by $\textnormal{\bf Scheme-3}$ and $P^2$, and $\textnormal{\bf Scheme-3F}$ and $P^3$ with upwind flux.  As predicted by the theorems, the total energy is not conserved exactly, but demonstrates good long time behavior. The fourth order scheme does a better job in conservation due to the higher order accuracy.

\begin{figure}[!htbp]
\centering
 \subfigure[Landau damping. ]{\includegraphics[width=0.45\textwidth]{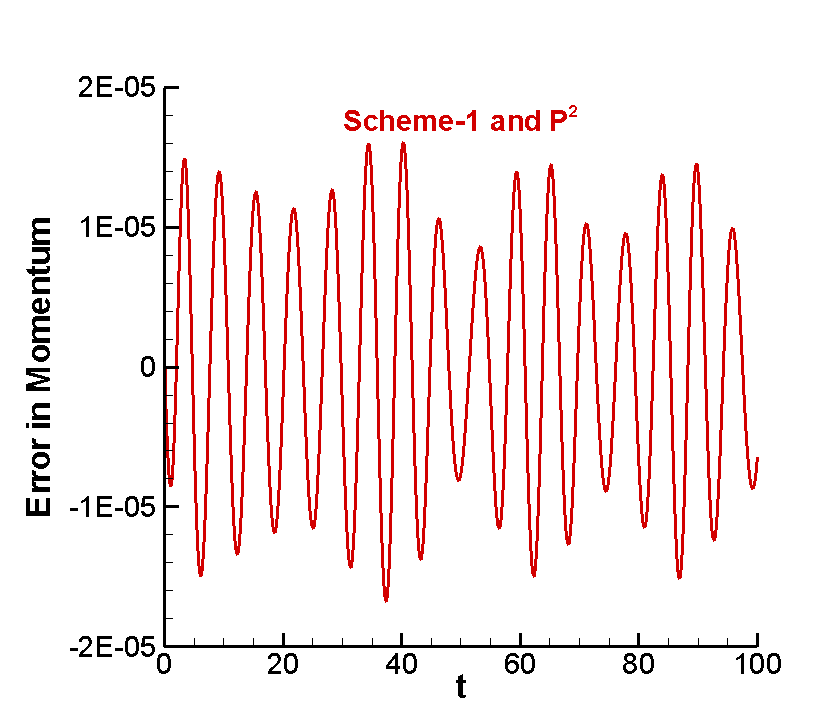}}
 \subfigure[Landau damping. ]{\includegraphics[width=0.45\textwidth]{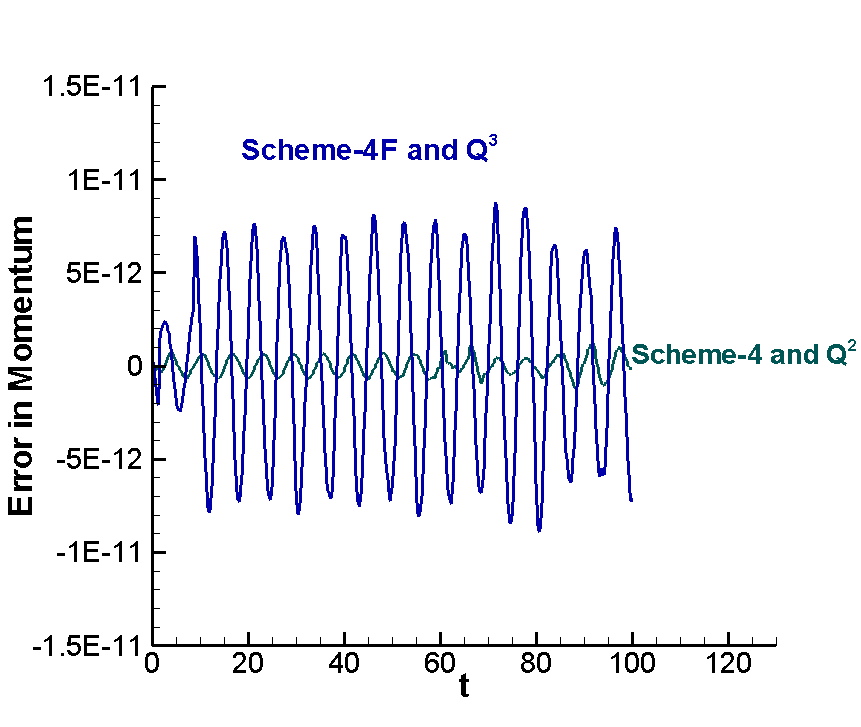}}\\
  \subfigure[Two-stream instability. ]{\includegraphics[width=0.45\textwidth]{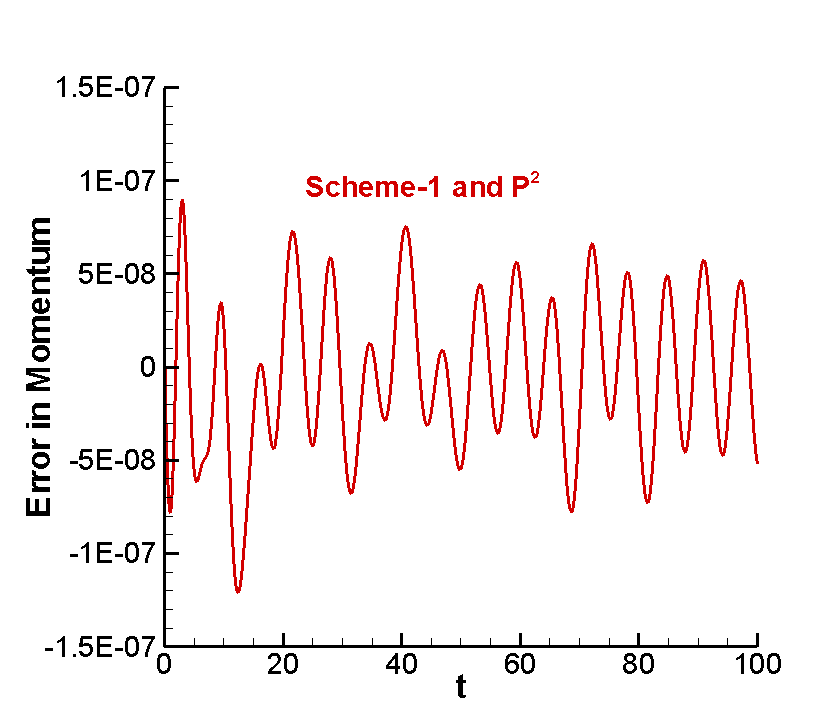}}
 \subfigure[Two-stream instability. ]{\includegraphics[width=0.45\textwidth]{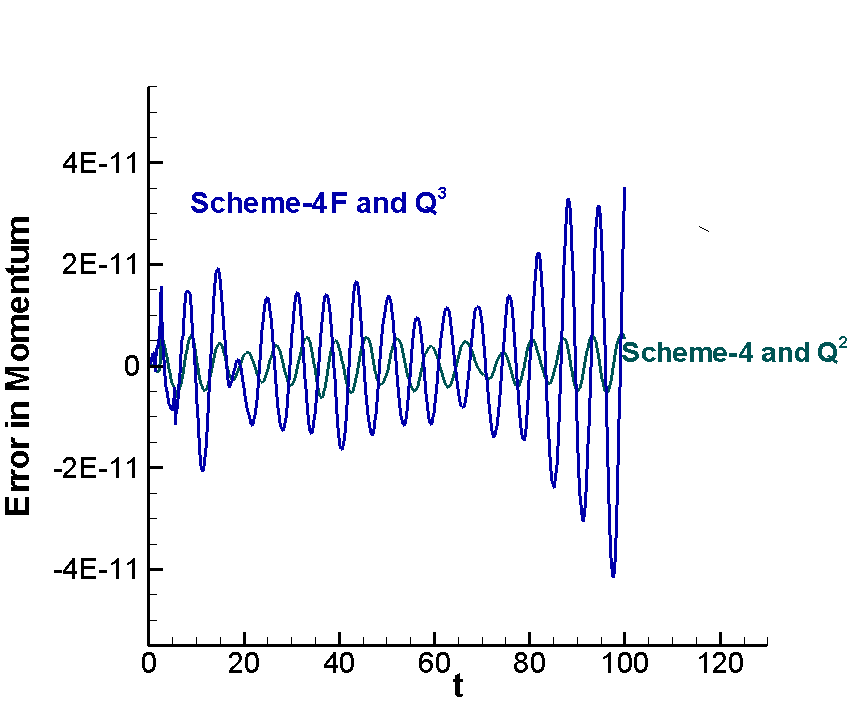}}
 \caption{Evolution of error in momentum by the indicated methods. $100\times200$ mesh. Upwind flux. }
 \label{figure_momentum}
\end{figure}

As for the example of Landau damping and two-stream instability, the momentum is also conserved.
In Figure \ref{figure_momentum}, we plot the relative error of momentum with the upwind flux. We can see the explicit scheme $\textnormal{\bf Scheme-1}$ gives relatively large errors for momentum, but the errors for  implicit schemes $\textnormal{\bf Scheme-4}$ and $\textnormal{\bf Scheme-4F}$ stay relatively small. While there is no rigorous proof, we believe this is due to the different treatment of the $f$ and $E$ component in $\textnormal{\bf Scheme-1}$.

\begin{figure}[!htbp]
\centering
 \subfigure[Landau damping. ]{\includegraphics[width=0.45\textwidth]{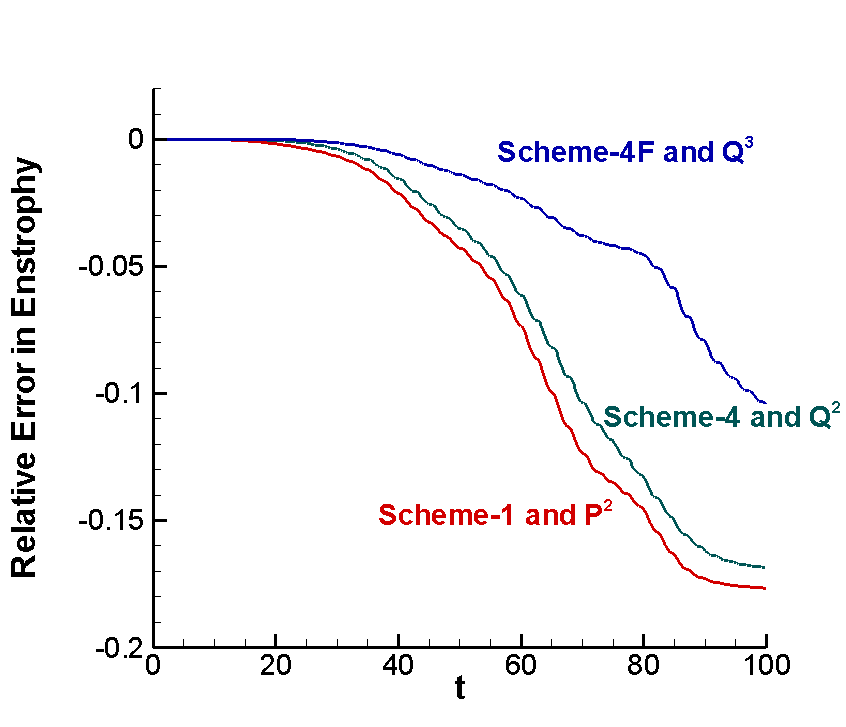}}
 \subfigure[Two-stream instability. ]{\includegraphics[width=0.45\textwidth]{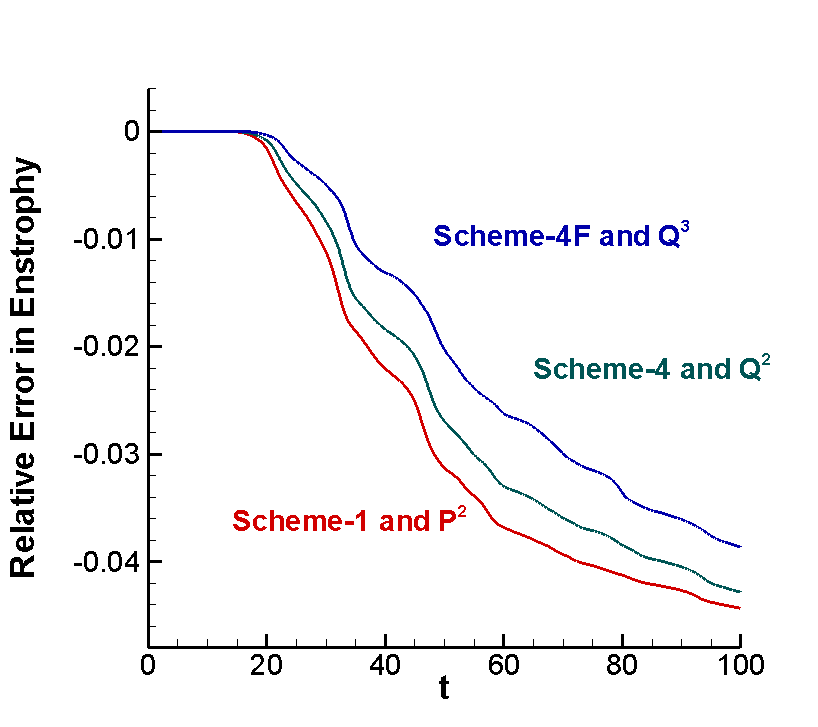}}\\
  \subfigure[Bump-on-tail instability. ]{\includegraphics[width=0.45\textwidth]{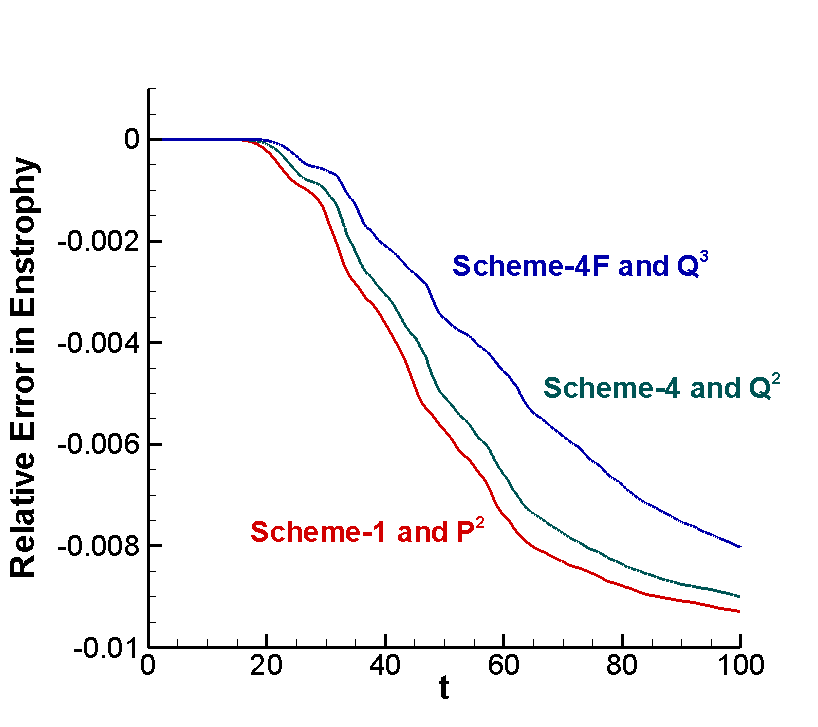}}
 \caption{Evolution of the relative error in enstrophy  by the indicated methods. $100\times200$ mesh. Upwind flux. }
 \label{figure_ens}
\end{figure}

In Figure \ref{figure_ens}, we plot the relative error of enstrophy. For upwind flux, the enstrophy is not conserved, and this is reflected in the figure. In particular, since  $\textnormal{\bf Scheme-4F}\textnormal{ and } Q^3$ is a fourth order scheme and more accurate than the other two methods, the relative error stays smaller compared to the other two lower order schemes.

 In Figure \ref{figure_mesh40}, we use a coarse mesh ($40 \times 80$) to plot the relative error in the conserved quantities to demonstrate that the conservation properties of our schemes are mesh independent. We use Landau damping, and $\textnormal{\bf Scheme-4} \textnormal{ and } Q^2$ to demonstrate the behavior. Upon comparison with the results from finer mesh  in Figures \ref{figure_energy1}, \ref{figure_momentum} and \ref{figure_ens}, we conclude that the mesh size has no impact on the conservation of total particle number and total energy as predicted by the theorems in Section \ref{sec:method}. This demonstrates the distinctive feature of our scheme: the total particle number and energy can be well preserved even with an under-resolved mesh.

\subsection{Comparison of central and upwind fluxes}
In this subsection, we use Landau damping as an example to compare the numerical solutions obtained by central and upwind fluxes. It is well known that for linear transport equations, the DG methods with upwind flux could achieve optimal $(k+1)$-th accuracy, but the central flux only obtain sub-optimal $k$-th order accuracy when $k$ is odd. Therefore, in our comparisons below, we only choose to use $k=2$ to compare the performance of the two fluxes.

\begin{figure}[!htbp]
\centering
 \subfigure[Total particle number.] {\includegraphics[width=0.45\textwidth]{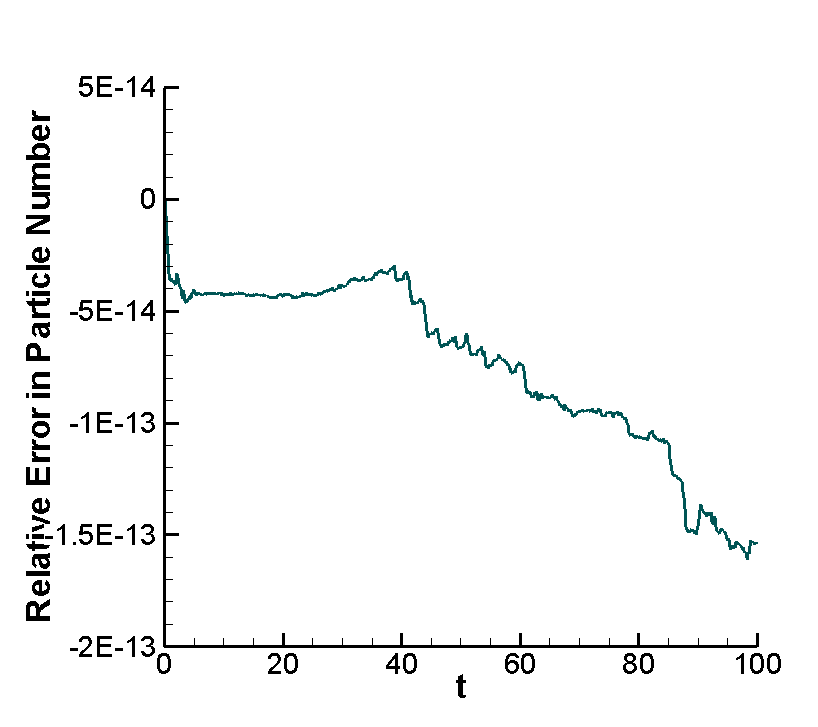}}
  \subfigure[Total energy. ]{\includegraphics[width=0.45\textwidth]{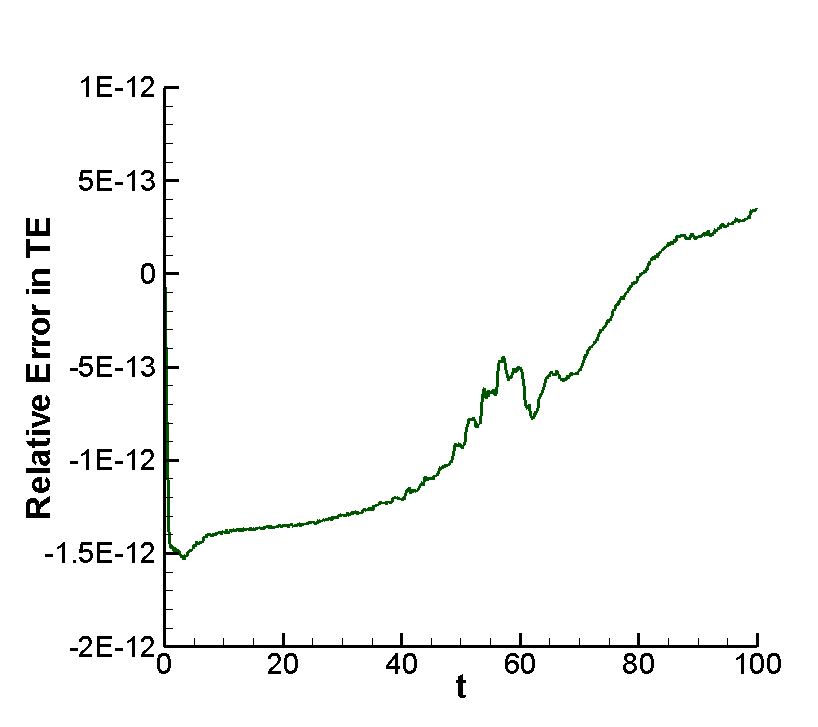}}\\
 \subfigure[Momentum. ]{\includegraphics[width=0.45\textwidth]{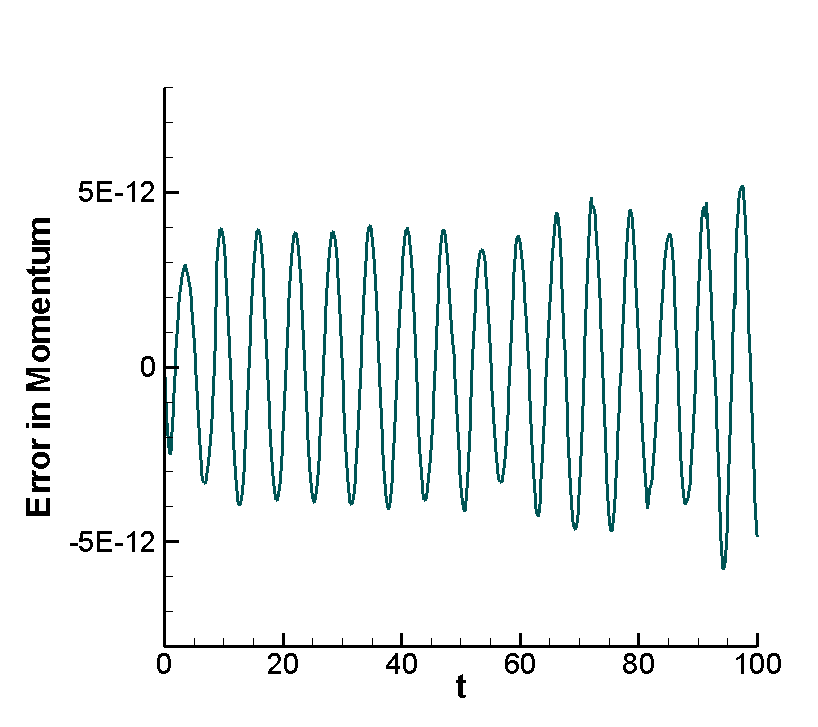}}
  \subfigure[Enstrophy. ]{\includegraphics[width=0.45\textwidth]{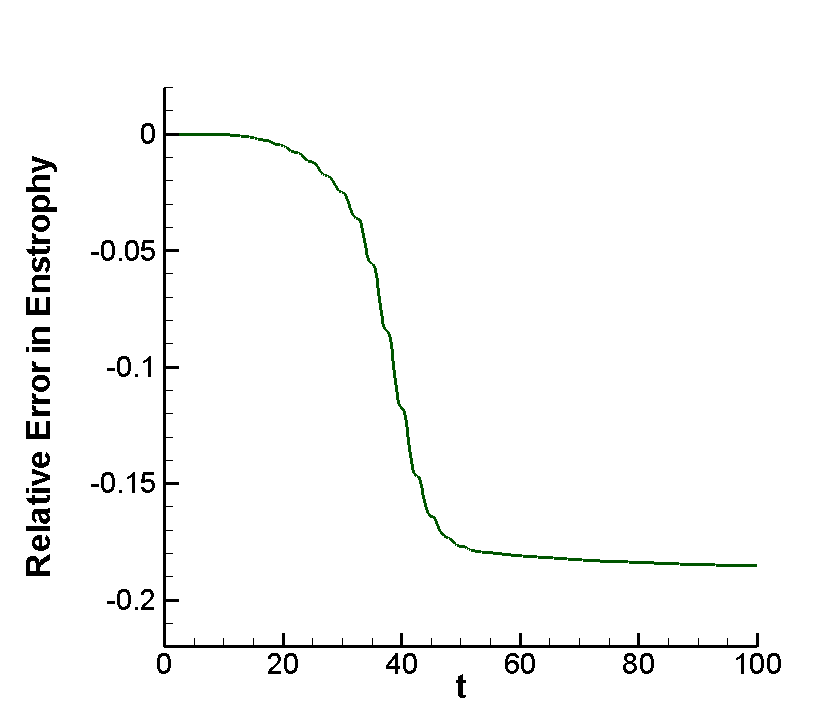}}
 \caption{Evolution of the relative error in conserved quantities by $\textnormal{\bf Scheme-4}\textnormal{ and } Q^2$ for Landau damping. $40\times80$ mesh. $CFL = 2$. Upwind flux. }
\label{figure_mesh40}
\end{figure}

\begin{figure}[!htbp]
\centering
 \subfigure[Total particle number.] {\includegraphics[width=0.45\textwidth]{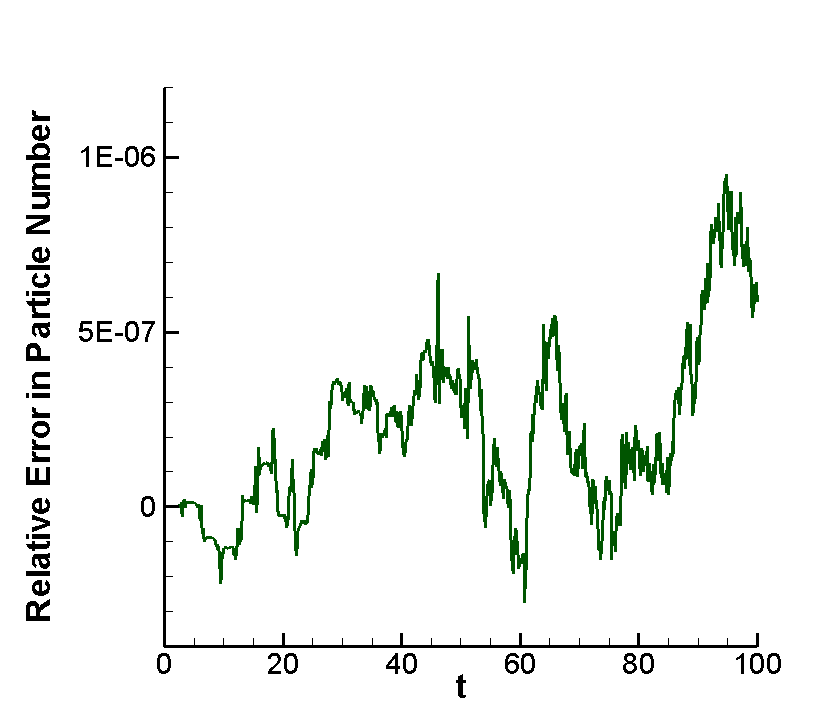}}
  \subfigure[Total energy. ]{\includegraphics[width=0.45\textwidth]{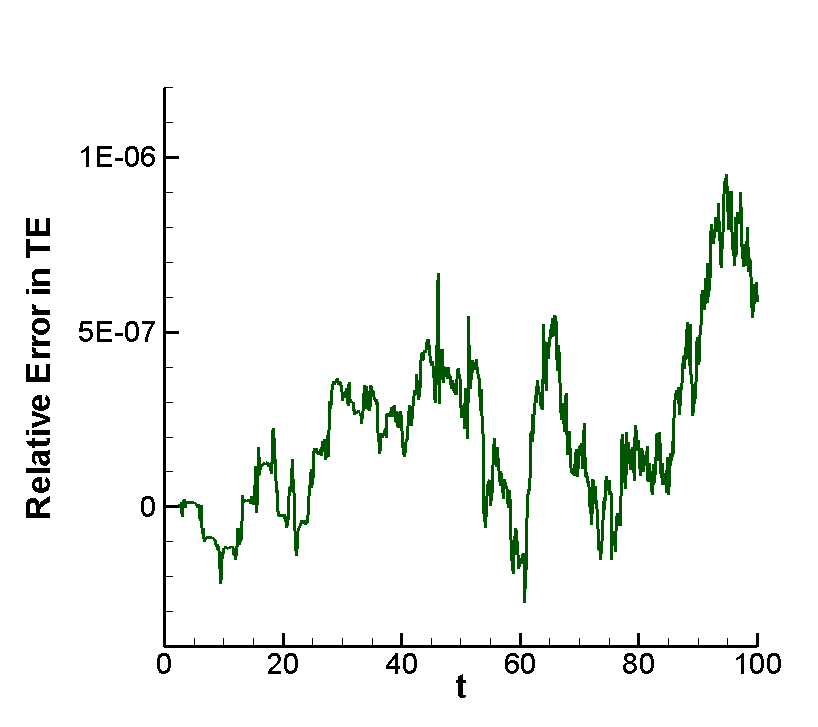}}\\
 \subfigure[Momentum. ]{\includegraphics[width=0.45\textwidth]{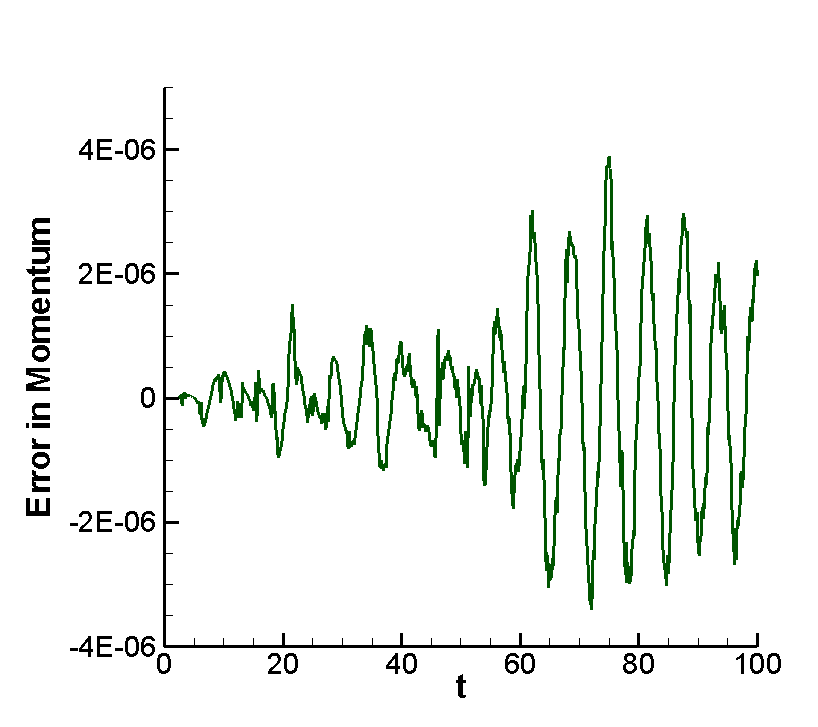}}
  \subfigure[Enstrophy. ]{\includegraphics[width=0.45\textwidth]{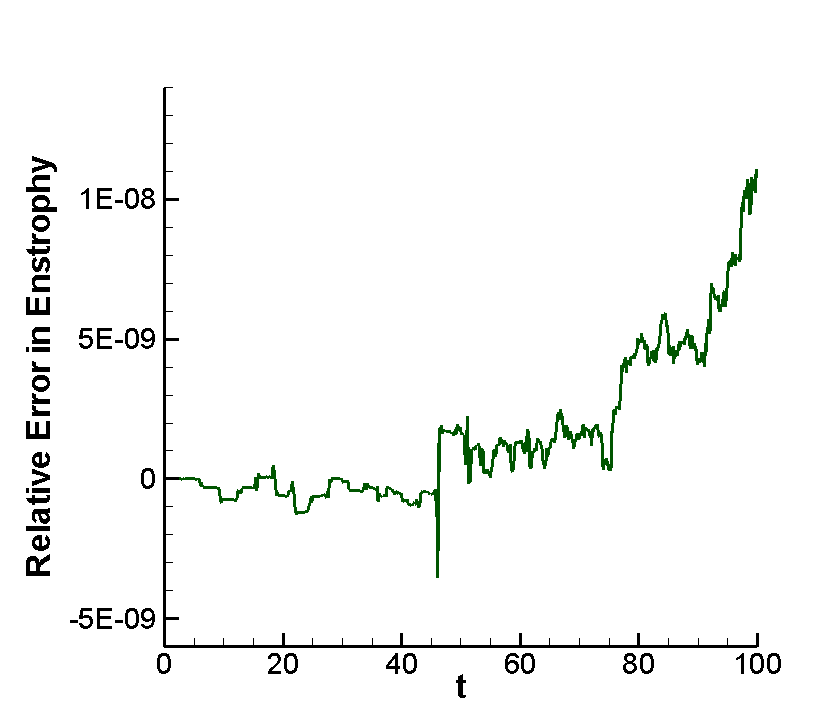}}
 \caption{Evolution of the relative error in conserved quantities by $\textnormal{\bf Scheme-4}\textnormal{ and } Q^2$ for Landau damping. $40\times80$ mesh. $CFL = 2$. Central flux. }
\label{figure_landaucentral}
\end{figure}

\begin{figure}[!htbp]
\centering
\subfigure[$t=10$. Upwind flux.]{\includegraphics[width=0.45\textwidth]{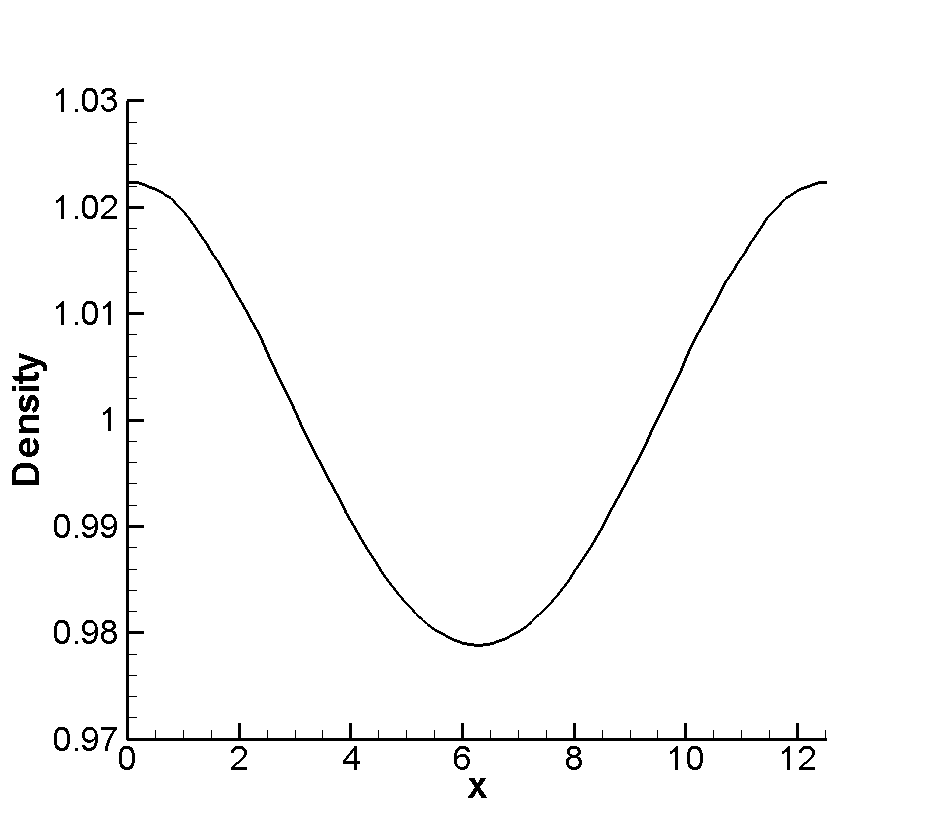}}
\subfigure[$t=10$. Central flux.]{\includegraphics[width=0.45\textwidth]{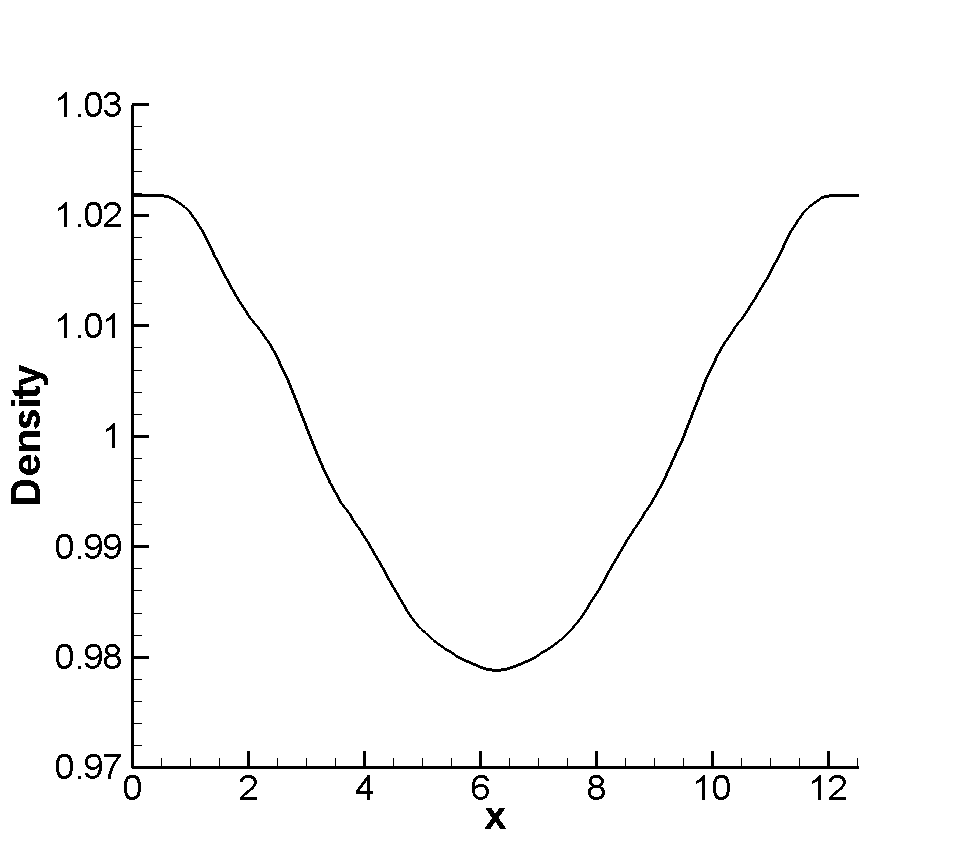}}\\
\subfigure[$t=50$. Upwind flux.]{\includegraphics[width=0.45\textwidth]{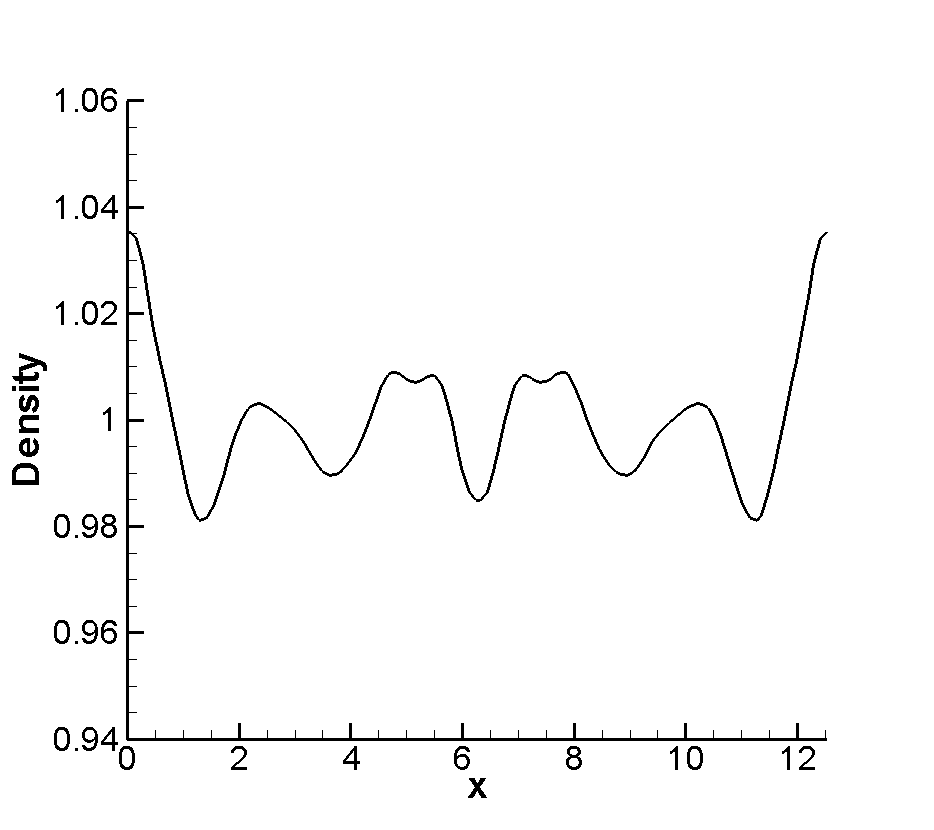}}
\subfigure[$t=50$. Central flux.]{\includegraphics[width=0.45\textwidth]{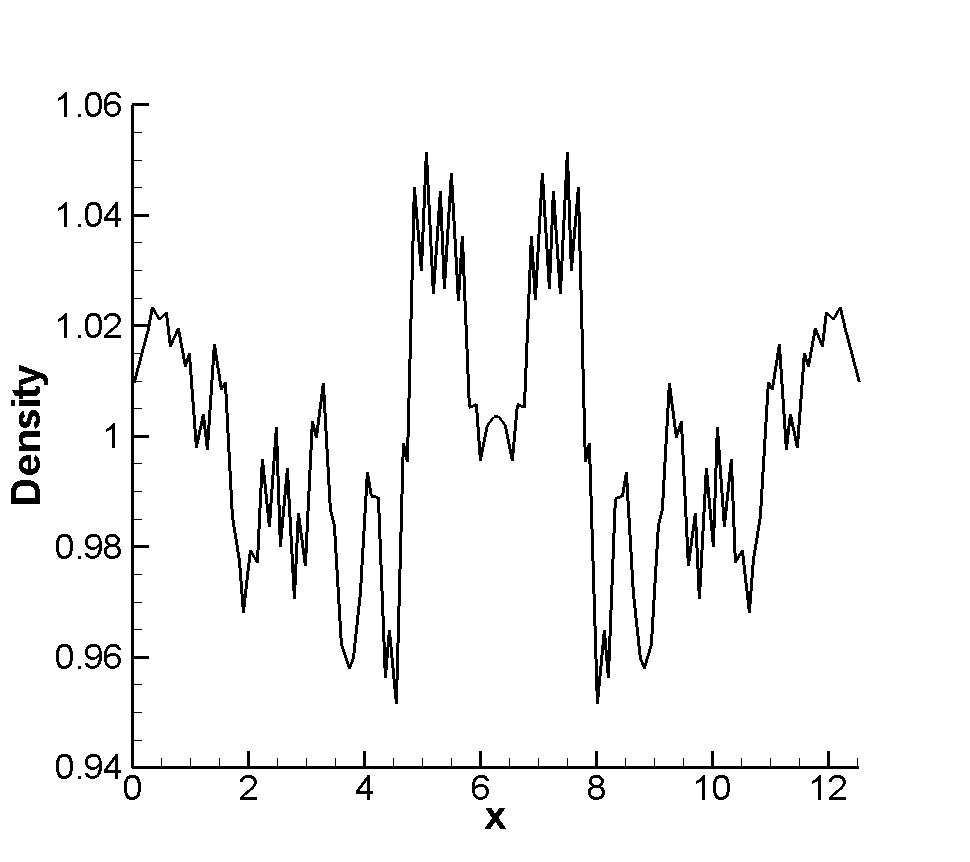}}\\
\subfigure[$t=90$. Upwind flux.]{\includegraphics[width=0.45\textwidth]{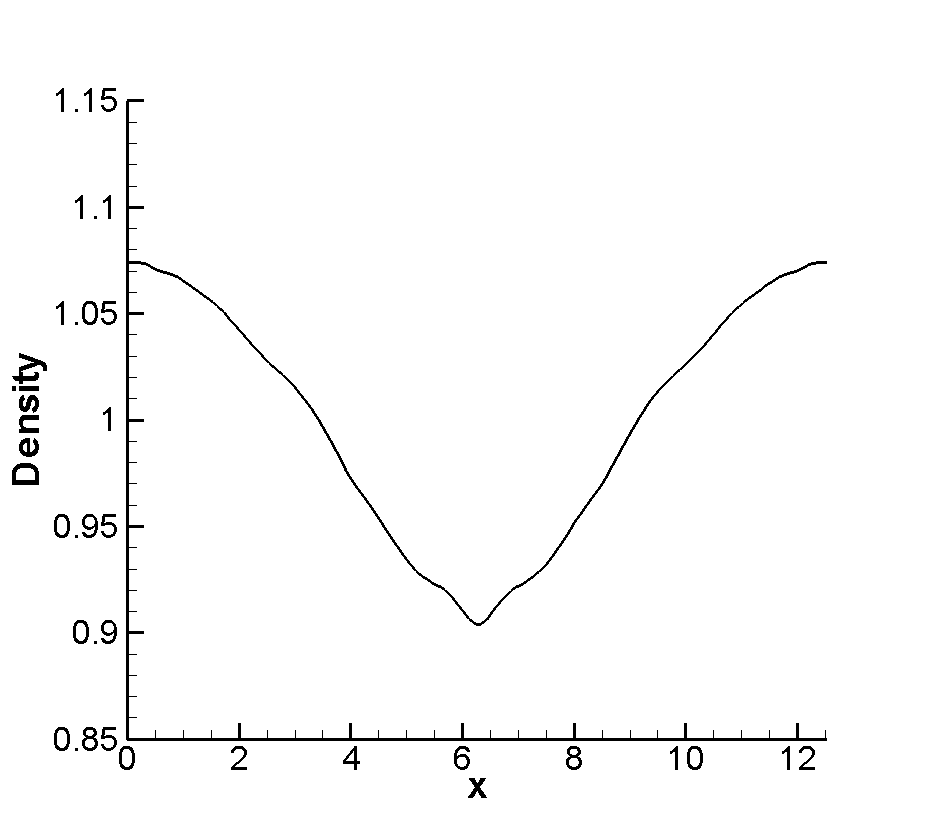}}
\subfigure[$t=90$. Central flux.]{\includegraphics[width=0.45\textwidth]{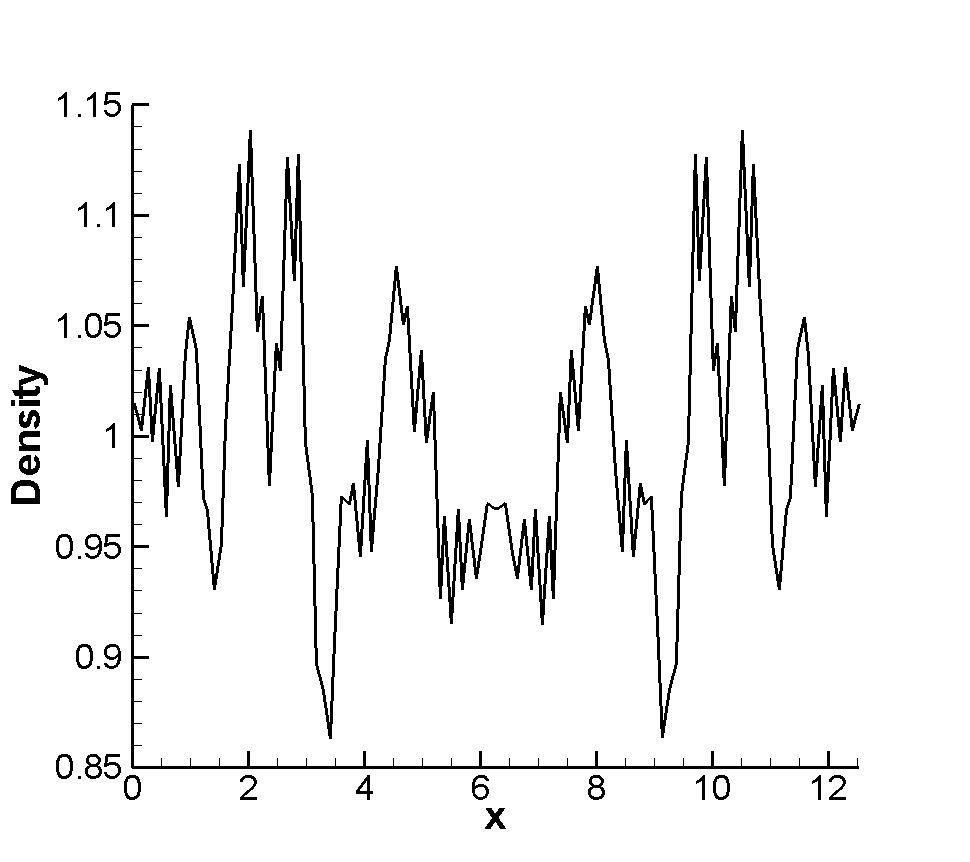}}\\
\caption{Comparison of density at the indicated time by upwind and central flux.  Landau damping. $\textnormal{\bf Scheme-4}\textnormal{ and } Q^2$. $40\times 80$ mesh. $CFL = 2$. }
\label{figure_landaudens}
\end{figure}

In Figure \ref{figure_landaucentral}, we list the numerical results obtained with scheme $\textnormal{\bf Scheme-4}\textnormal{ and } Q^2$, on a $40 \times 80$ mesh, with central flux. Upon comparison with
 Figure \ref{figure_mesh40}, we can see that the central flux could preserve the enstrophy much better, on the scale of $10^{-8}$. This is predicted in Section \ref{sec:method}. However, the conservation of particle number, momentum, and total energy are less satisfactory, about 6 magnitudes bigger than the upwind flux on the same mesh. The reason is because central flux does not build any numerical dissipation into the scheme, and this is not desired  when filamentation occurs. Lack of dissipation could make the numerical solution   high oscillatory, resulting in a non-zero value near the velocity boundaries, and causing the loss of conservation. This fact is illustrated further in Figure \ref{figure_landaudens}, where we observe the density generated by central flux is more oscillatory than the ones generated by the upwind flux.

\subsection{Convergence of the Newton-Krylov solver}
In this subsection, we investigate the relation of convergence of the Newton-Krylov solver and the $CFL$ number.  In particular, 
we implement scheme $\textnormal{\bf Scheme-4}$ on a $100\times200$ mesh and integrate up to $ T=50$ with polynomial space $\mS_h^1$, and set the tolerance parameter to be $\epsilon_{tol}=10^{-12}$. In Table \ref{CFL_landau}, $nni$ is the average of the numbers of  Newton iterations per step, and $nli$ is the average of the  number of Krylov iterations per step. We can see when the $CFL$ number increases, the number of Newton and Krylov iterations increase as expected; however, the increase seems to be sub-linear. This indicates that it is more efficient to use a larger $CFL$ when accuracy permits. We remark that the Newton-Krylov solver fails to converge if we increase the  $CFL$ to  $300$ in this case.
\begin{table}[!htbp]
\centering
\caption{Relation of the number of Newton and Krylov iterations per step with the $CFL$ number. $\textnormal{\bf Scheme-4}\textnormal{ and } Q^1$. $100\times 200$ mesh. $\epsilon_{tol}=10^{-12}$.}
\label{CFL_landau}
\medskip
\begin{tabular}{|c|c|c|c|c|c|c|c|c|c|}\hline
$CFL$	&1	&10	&20	&40&	80	&100	&150	&200	&250\\
\hline
$nni$	
&4.29
&4.94
&5.16
&5.50
&6.10
&6.79
&7.81
&8.97
&9.86\\
\hline
$nli$	
&6.26
&
13.36
&
19.82
&33.02
&59.36
&82.01
&106.64
&122.27
&147.71\\
\hline
\end{tabular}
\end{table}

\subsection{Collections of numerical data}
\begin{figure}[!htbp]
\centering
 \subfigure[log $FM_1$. ]{\includegraphics[width=0.45\textwidth]{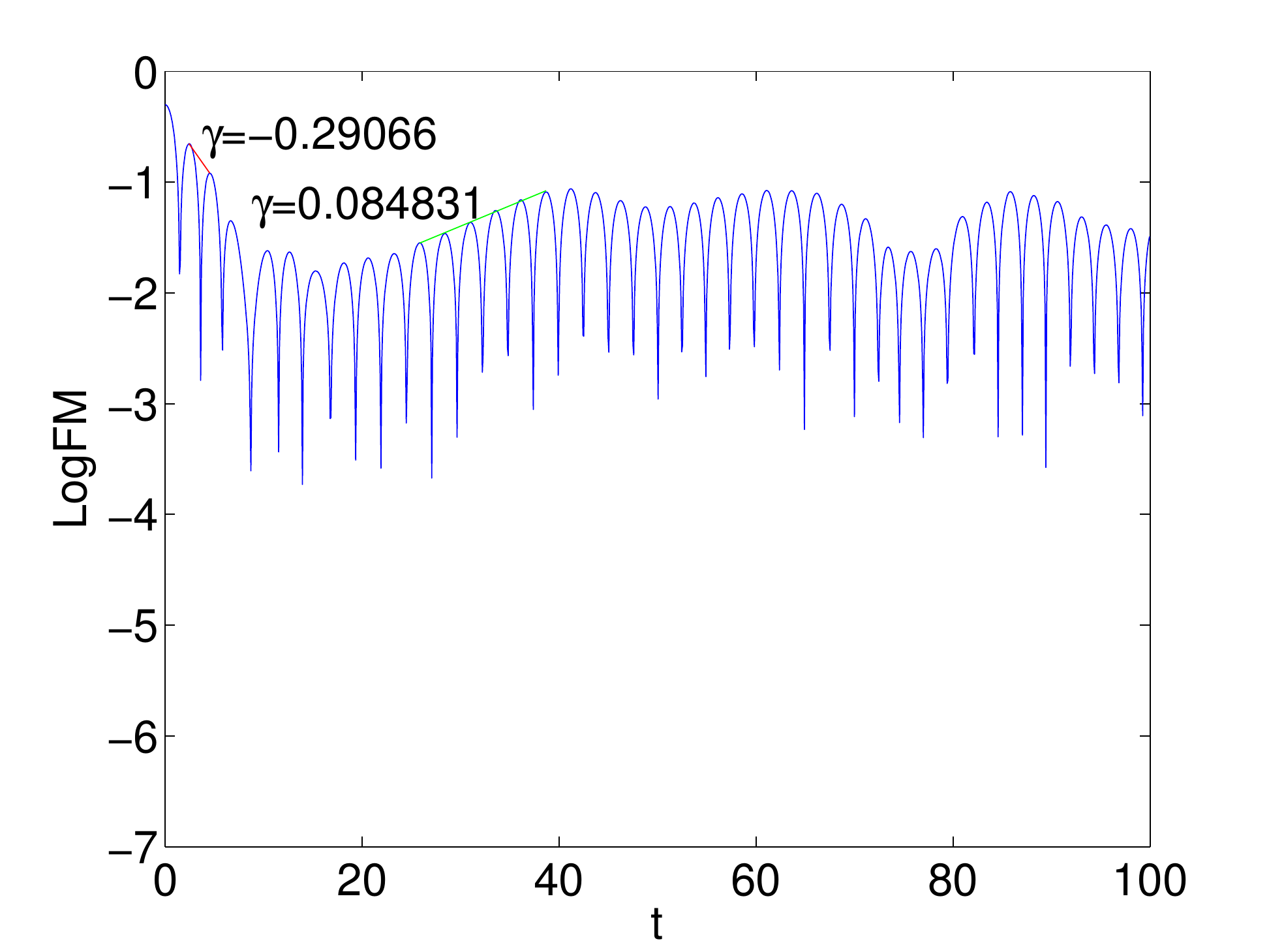}}
 \subfigure[log $FM_2$. ]{\includegraphics[width=0.45\textwidth]{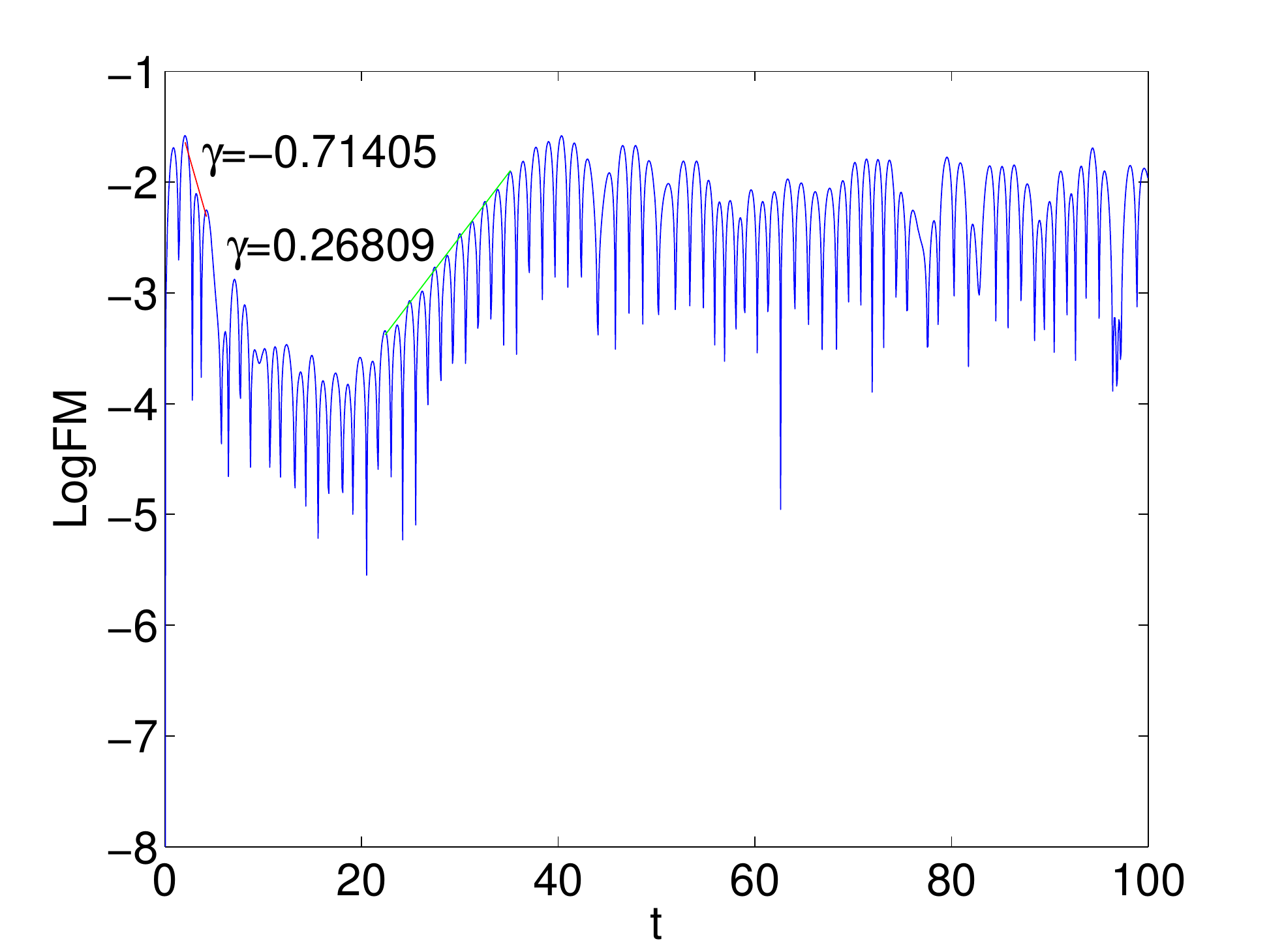}}\\
  \subfigure[log $FM_3$.  ]{\includegraphics[width=0.45\textwidth]{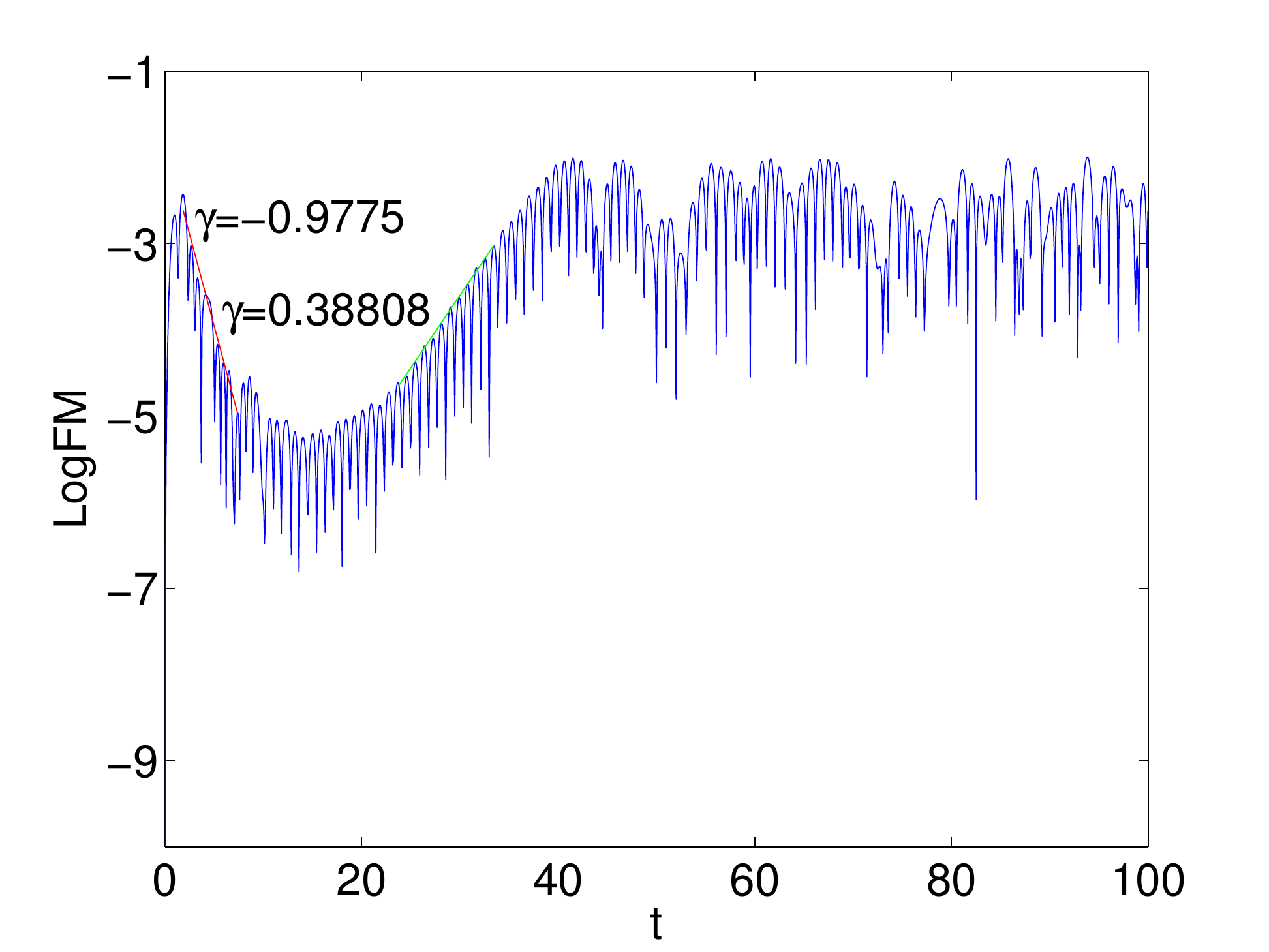}}
  \subfigure[log $FM_4$.]{\includegraphics[width=0.45\textwidth]{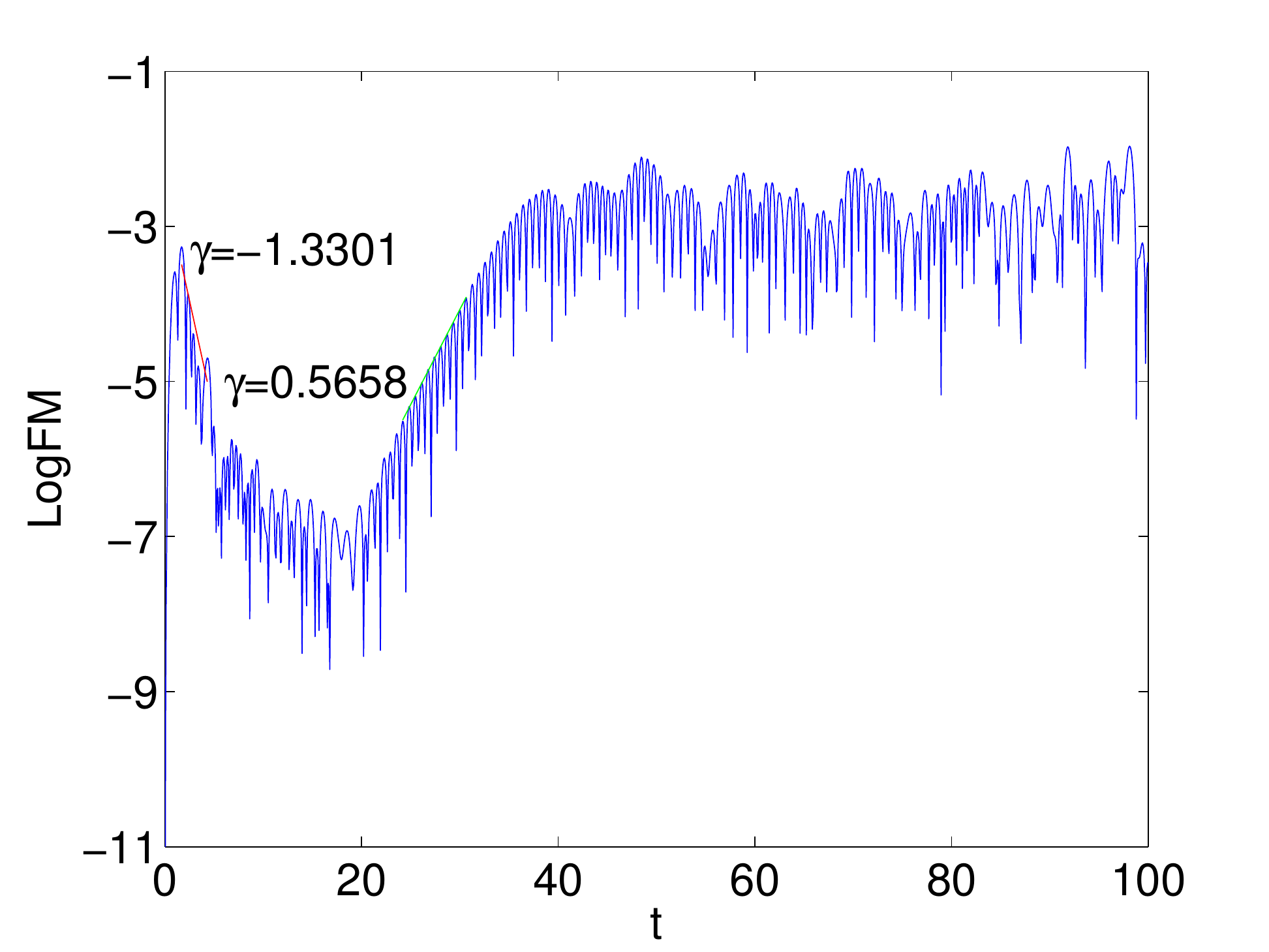}}
 \caption{Log Fourier modes   of Landau damping. $\textnormal{\bf Scheme-4F}\textnormal{ and } Q^3$. $100\times200$ mesh. Upwind flux. }
\label{figure_landaufm}
\end{figure}

\begin{figure}[!htbp]
\centering
\subfigure[Two-stream instability.]{\includegraphics[width=0.45\textwidth]{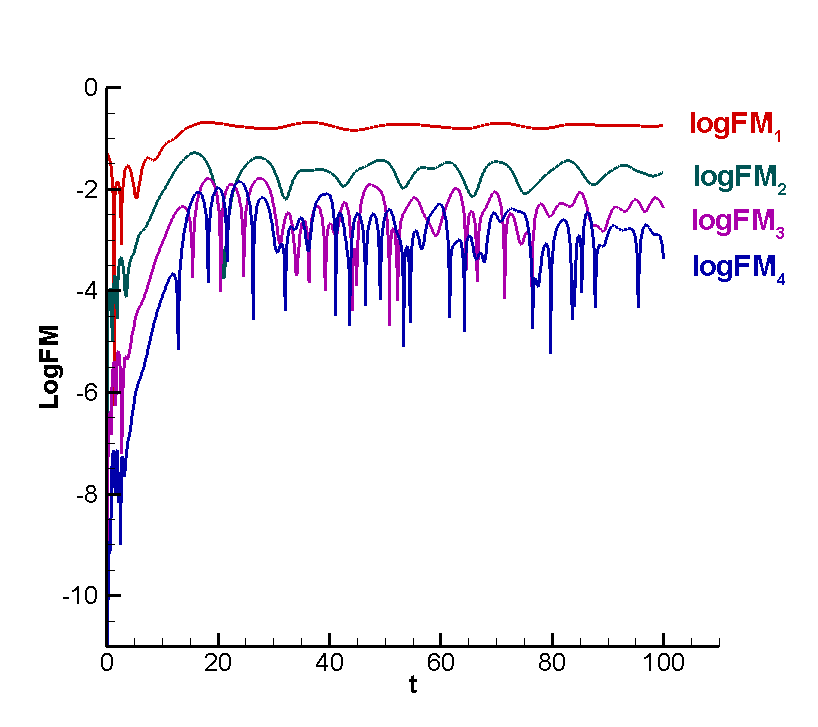}}
\subfigure[Bump-on-tail instability.]{\includegraphics[width=0.45\textwidth]{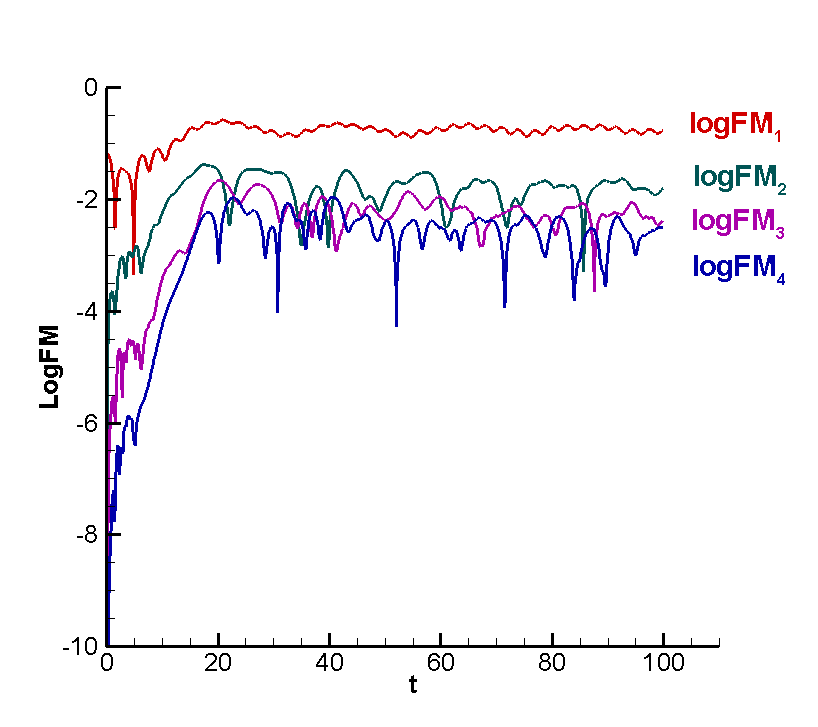}}
\caption{Log Fourier modes. $\textnormal{\bf Scheme-4F}\textnormal{ and } Q^3$. $100\times 200$ mesh. Upwind flux.}
\label{figure_tsbumpfm}
\end{figure}

\begin{figure}[!htbp]
\centering
\subfigure[$t=0$.]{\includegraphics[width=0.45\textwidth]{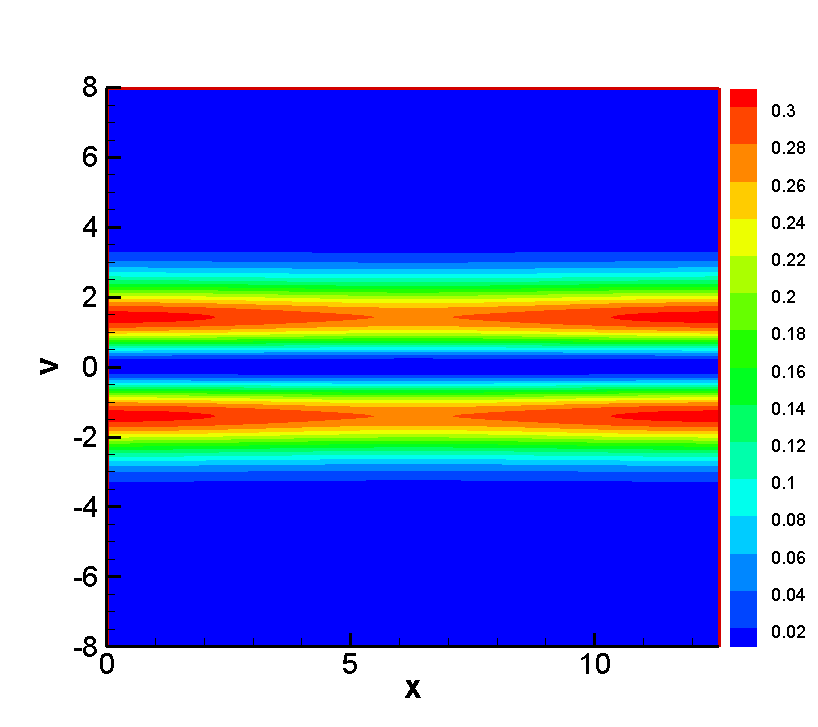}}
\subfigure[$t=10$.]{\includegraphics[width=0.45\textwidth]{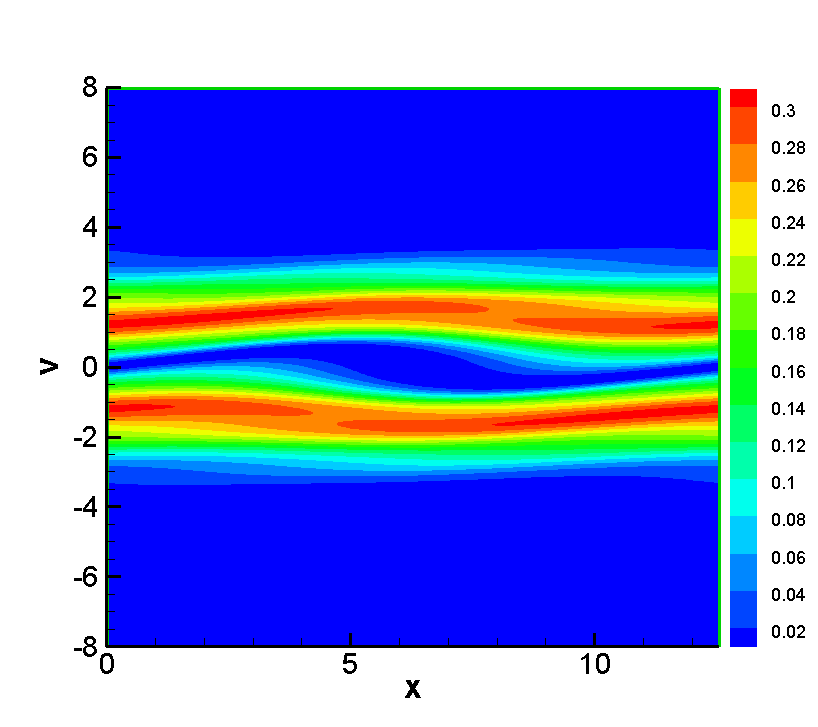}}\\
\subfigure[$t=20$.]{\includegraphics[width=0.45\textwidth]{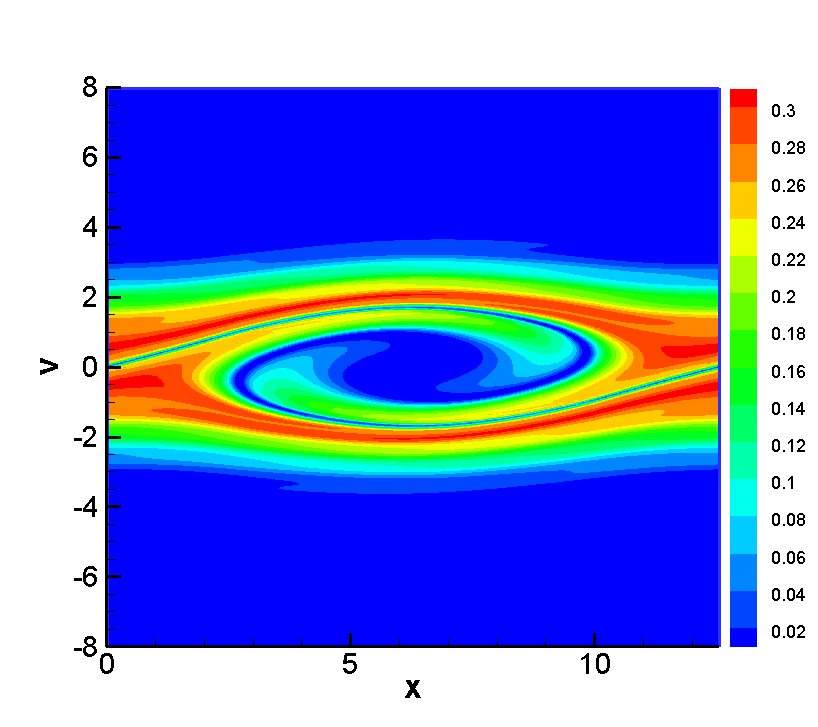}}
\subfigure[$t=40$.]{\includegraphics[width=0.45\textwidth]{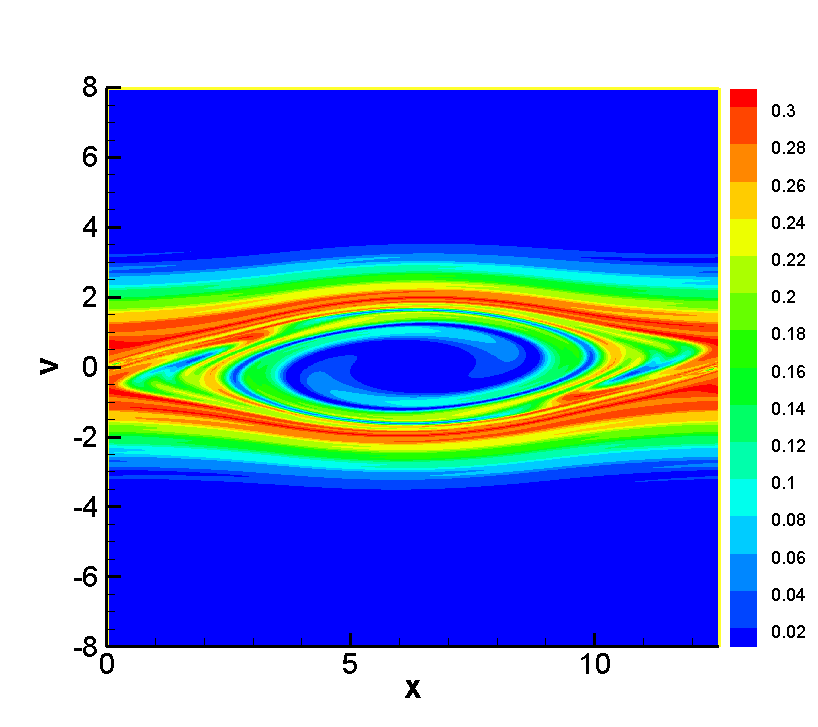}}\\
\subfigure[$t=80$.]{\includegraphics[width=0.45\textwidth]{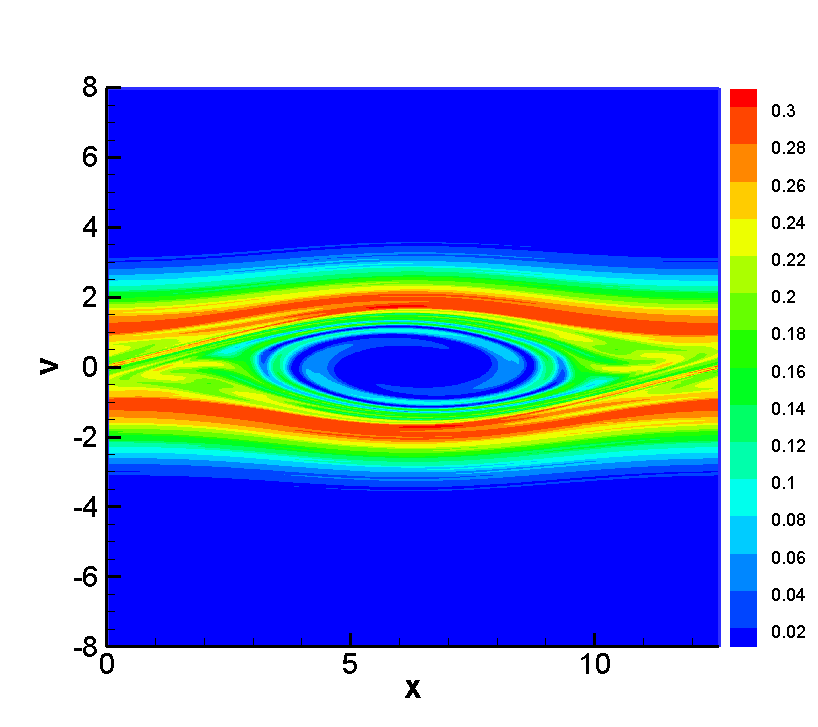}}
\subfigure[$t=100$.]{\includegraphics[width=0.45\textwidth]{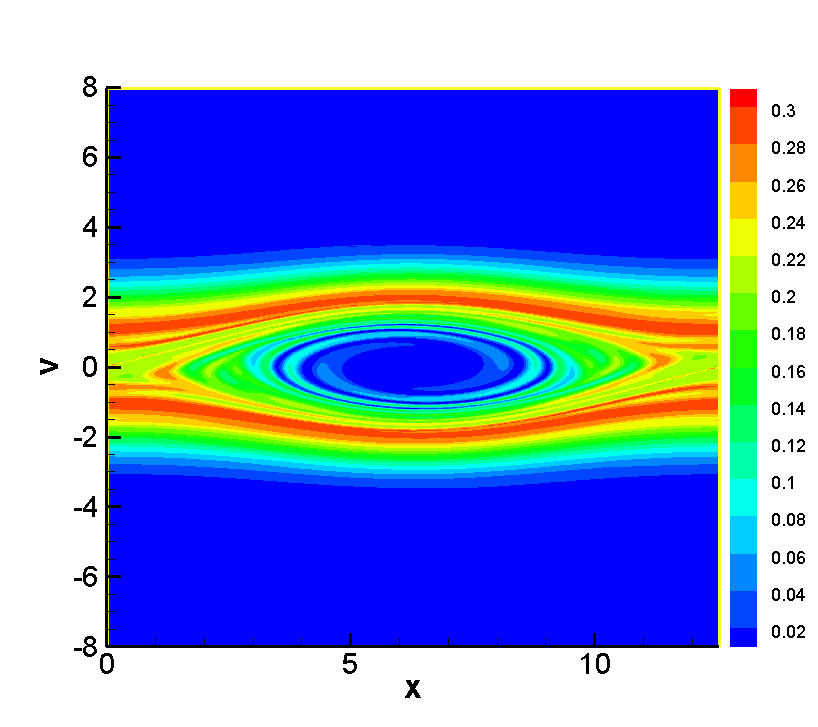}}
\caption{Phase space contour plots for the two-stream instability at the   indicated time. $\textnormal{\bf Scheme-4F}\textnormal{ and } Q^3$. $100\times 200$ mesh. Upwind flux.}
\label{figure_tscontour}
\end{figure}

\begin{figure}[!htbp]
\centering
\subfigure[$t=0$.]{\includegraphics[width=0.45\textwidth]{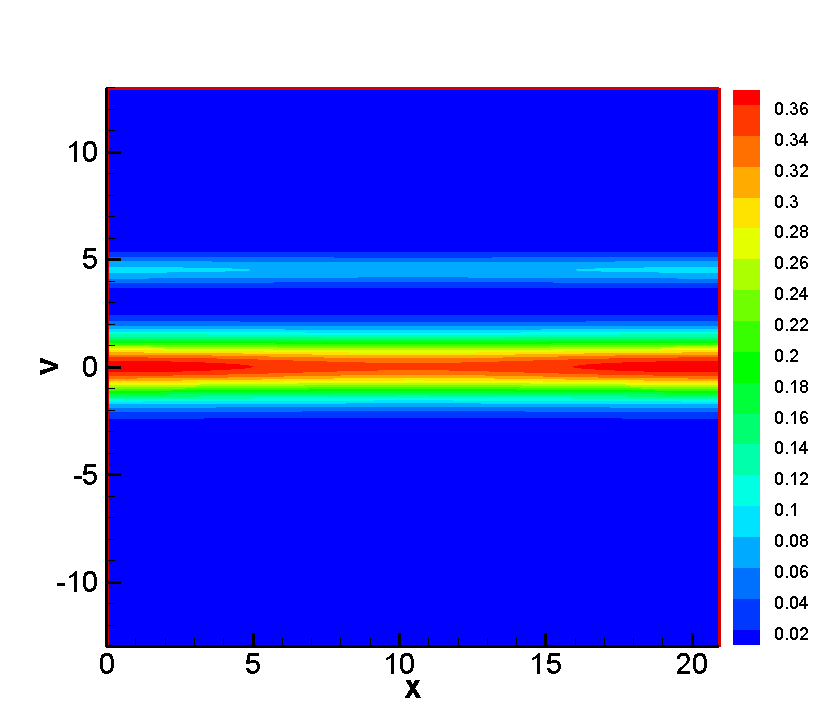}}
\subfigure[$t=15$.]{\includegraphics[width=0.45\textwidth]{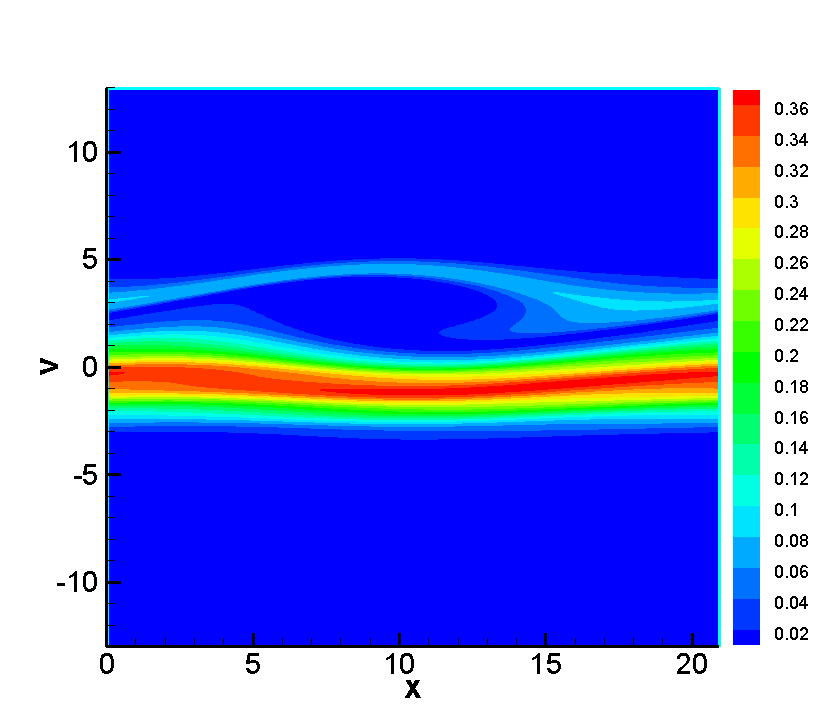}}\\
\subfigure[$t=20$.]{\includegraphics[width=0.45\textwidth]{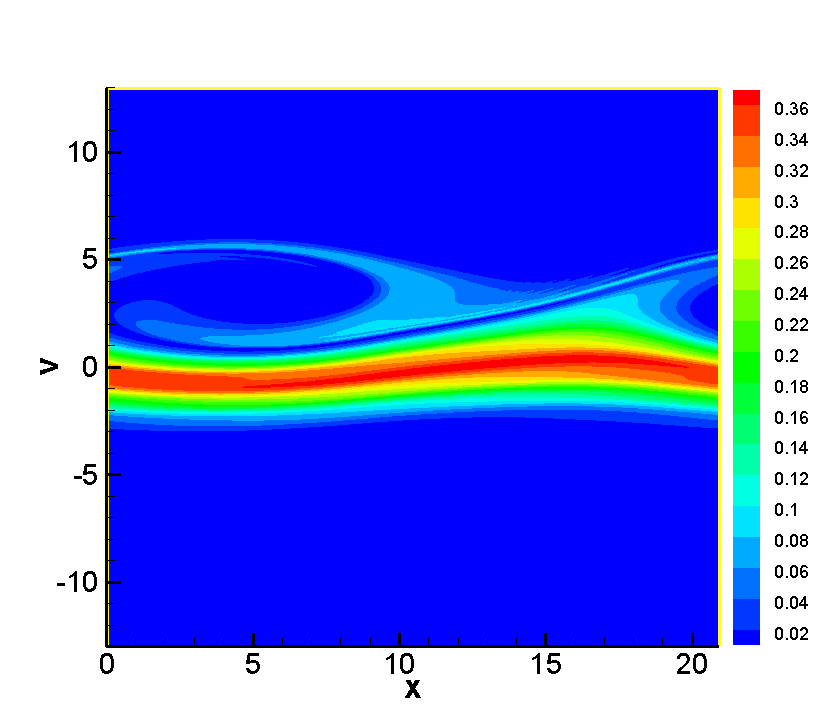}}
\subfigure[$t=30$.]{\includegraphics[width=0.45\textwidth]{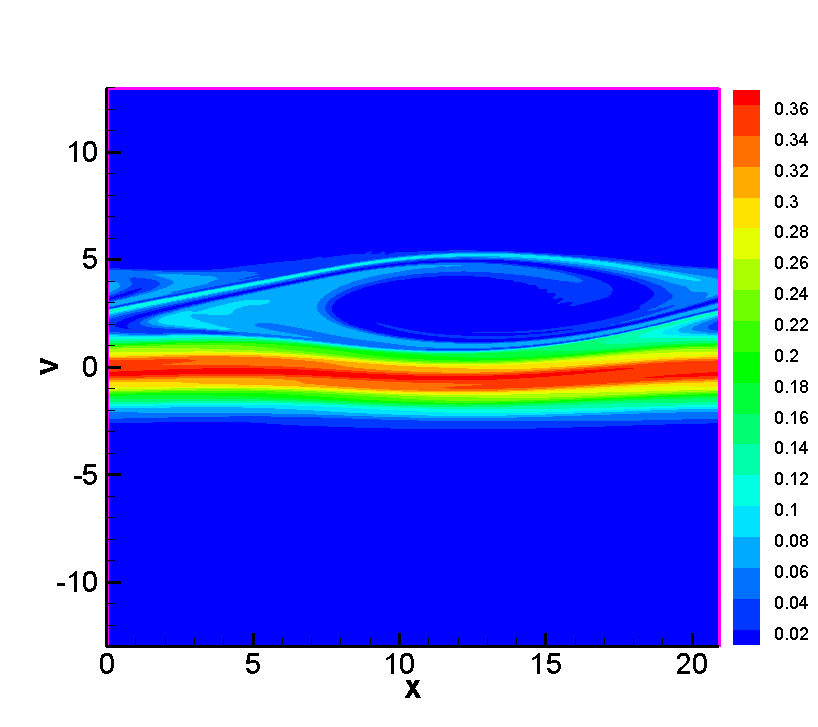}}\\
\subfigure[$t=50$.]{\includegraphics[width=0.45\textwidth]{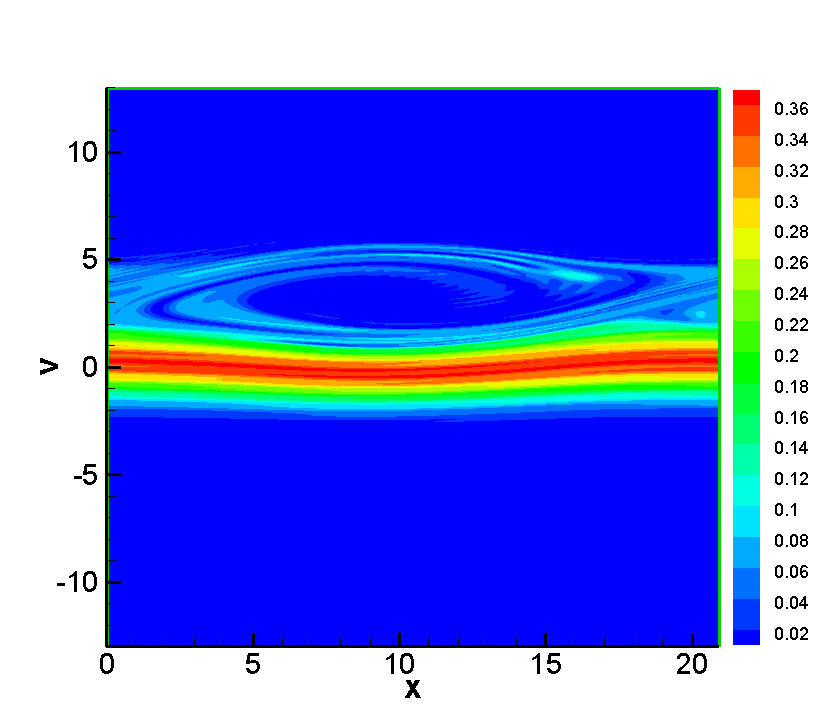}}
\subfigure[$t=60$.]{\includegraphics[width=0.45\textwidth]{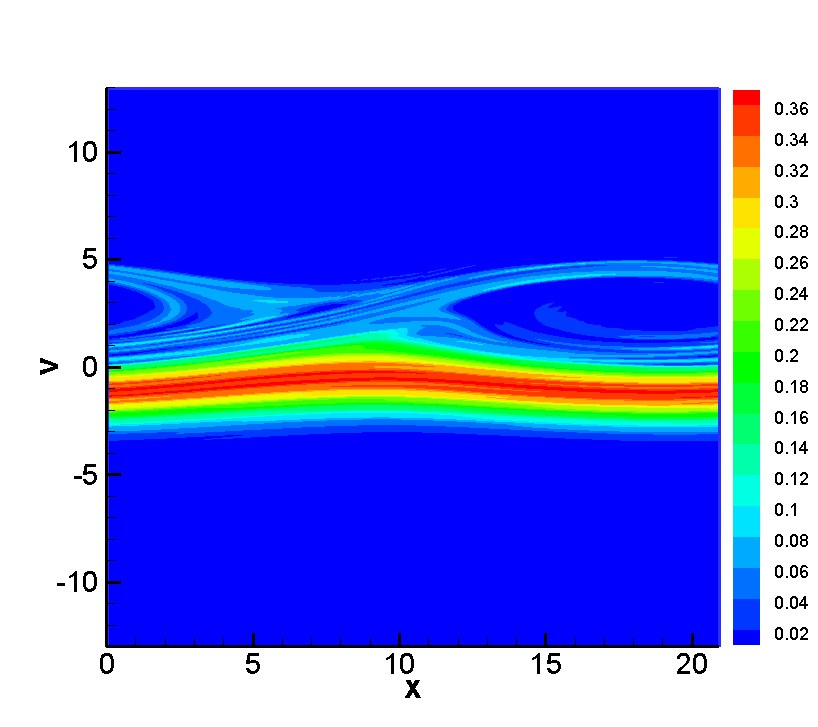}}
\caption{Phase space contour plots for the bump-on-tail instability at the  indicated time. $\textnormal{\bf Scheme-4F}\textnormal{ and } Q^3$. $100\times 200$ mesh. Upwind flux.}
\label{figure_bumpcontour}
\end{figure}

\begin{figure}[!htbp]
\centering
\subfigure[$t=0$.]{\includegraphics[width=0.45\textwidth]{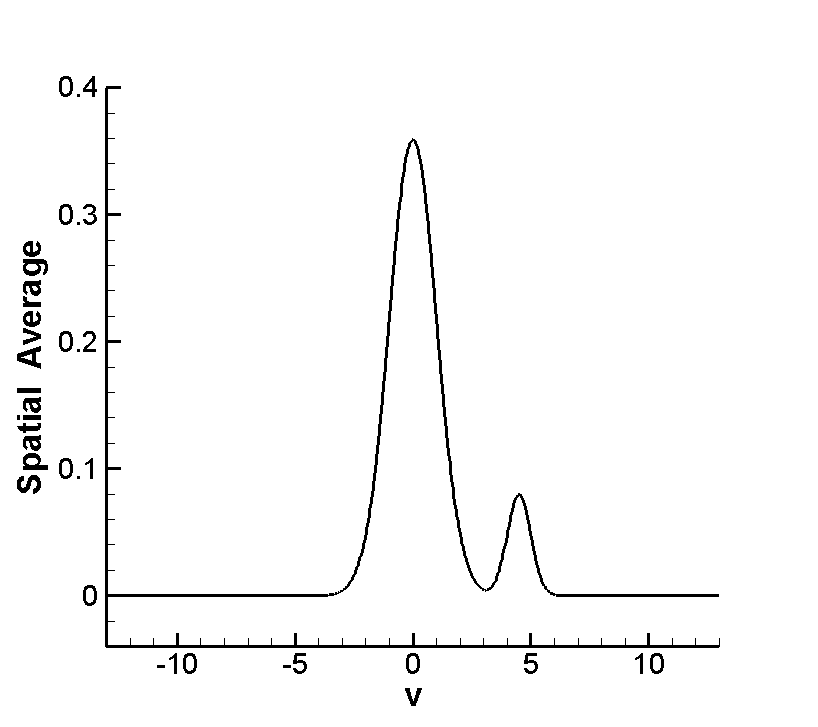}}
\subfigure[$t=15$.]{\includegraphics[width=0.45\textwidth]{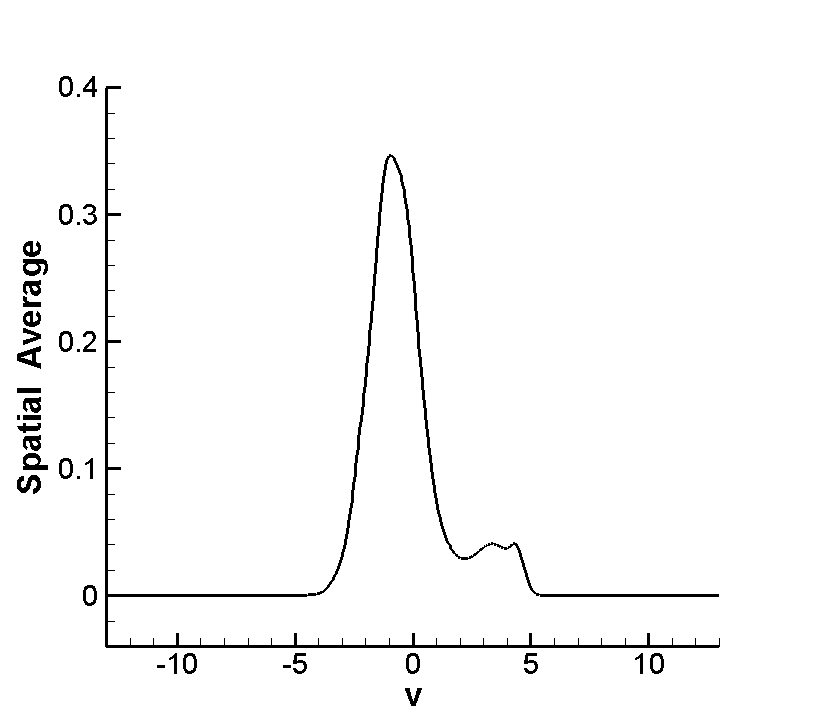}}\\
\subfigure[$t=20$.]{\includegraphics[width=0.45\textwidth]{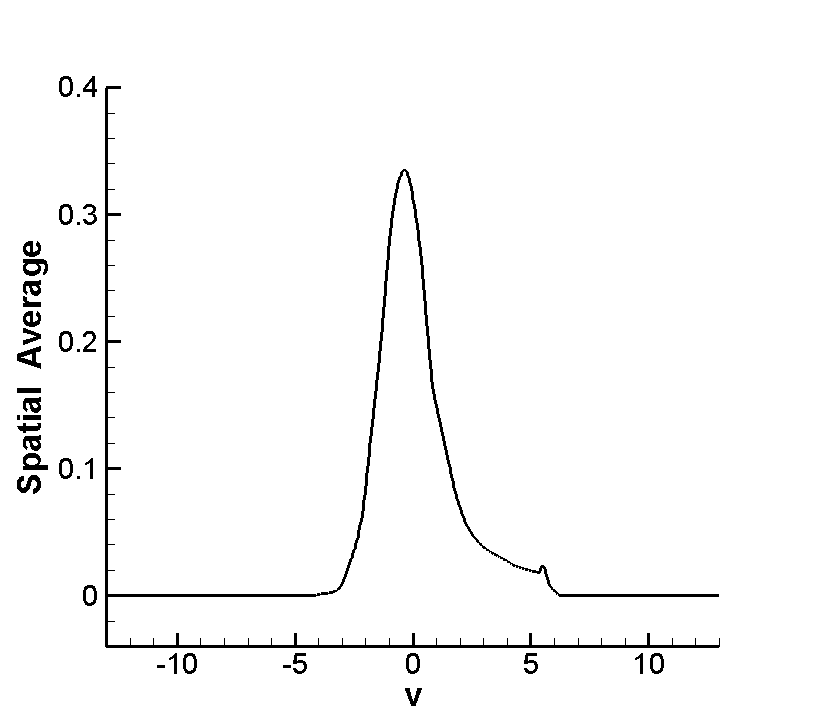}}
\subfigure[$t=30$.]{\includegraphics[width=0.45\textwidth]{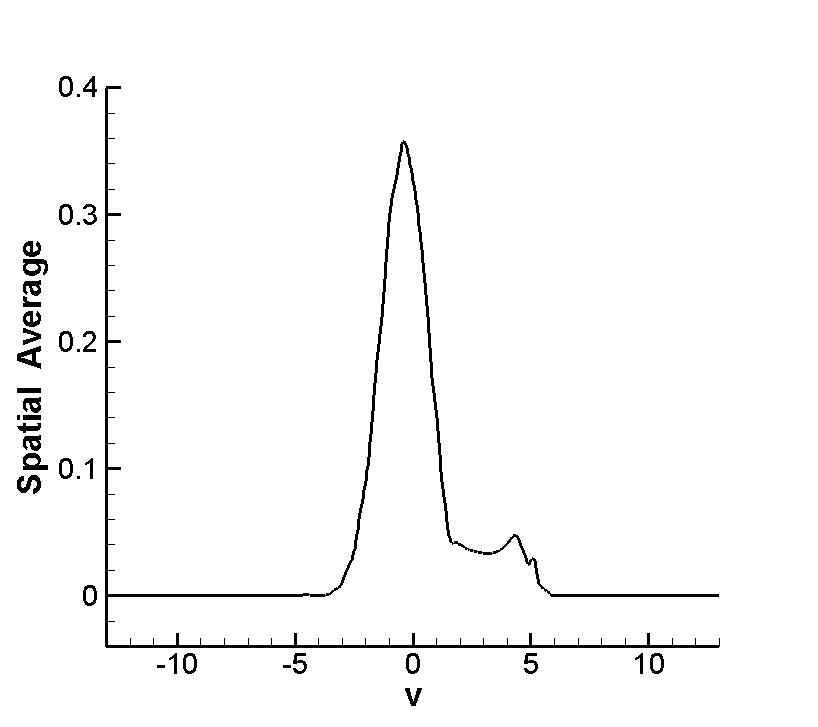}}\\
%
\subfigure[$t=50$.]{\includegraphics[width=0.45\textwidth]{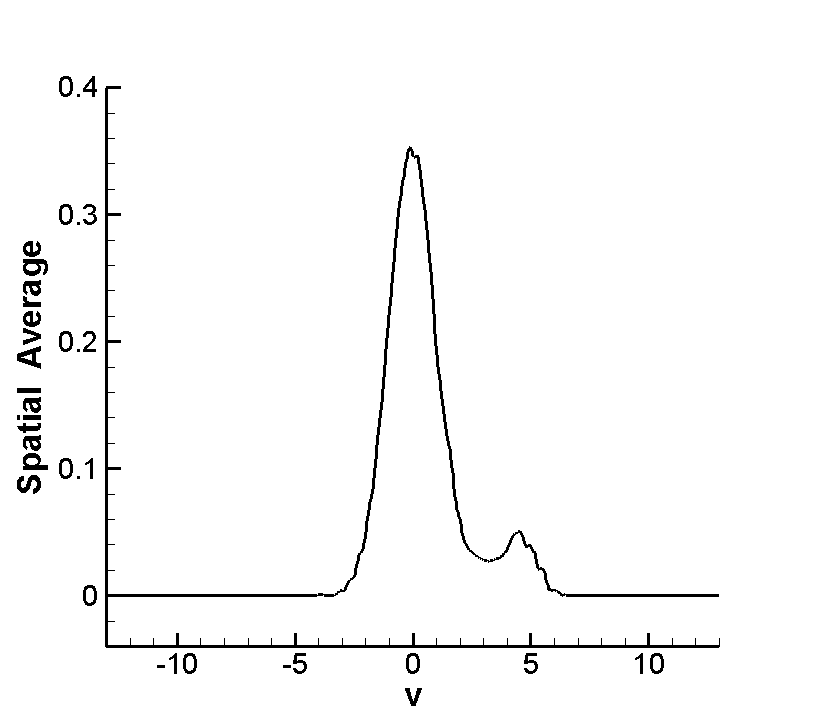}}
\subfigure[$t=60$.]{\includegraphics[width=0.45\textwidth]{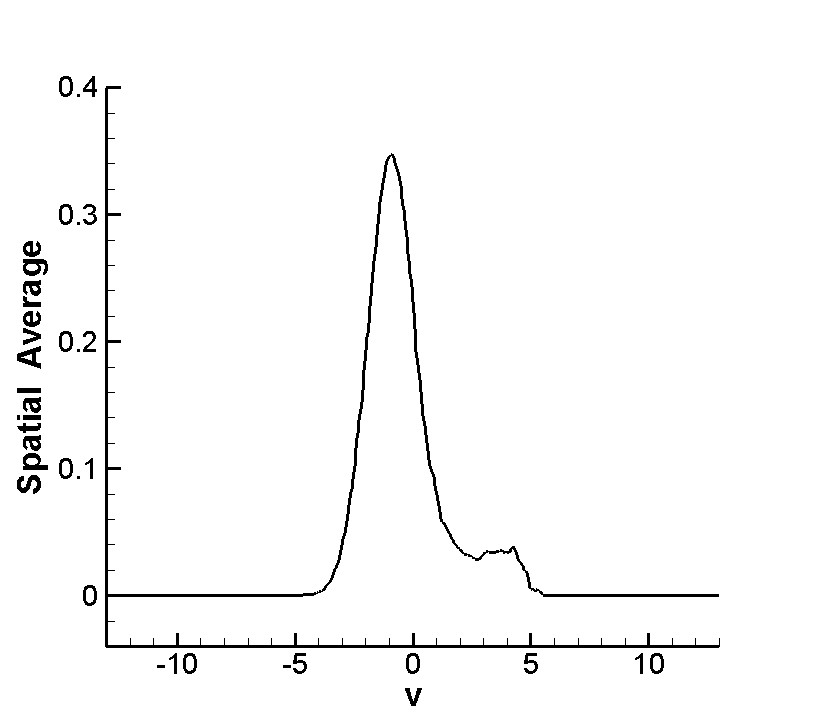}}
\caption{Spatial average of the distribution function for the bump-on-tail instability at the indicated time. $\textnormal{\bf Scheme-4F}\textnormal{ and } Q^3$. $100\times 200$ mesh. Upwind flux.}
\label{figure_bumpspa}
\end{figure}

In this subsection, we collect some sample numerical data to benchmark our schemes. In particular, we use the fourth order accurate scheme $\textnormal{\bf Scheme-4F}\textnormal{ and } Q^3$ on a $100\times200$ mesh with upwind flux.
In Figures \ref{figure_landaufm}, \ref{figure_tsbumpfm}, we plot the Log Fourier modes for the electric field for all three examples, where
the $n$-th Log Fourier mode for the electric field $E(x,t)$ \cite{Heath} is defined as
$$
logF\!M_n(t)=\log_{10} \left(\frac{1}{L} \sqrt{\left|\int_0^L E(x, t) \sin(\kappa nx) \, dx \right|^2 +
\left|\int_0^L E(x, t) \cos(\kappa nx)\,  dx \right|^2} \right).
$$
The Log Fourier modes generated by our methods agree well with other solvers in the literature.  In Figures \ref{figure_tscontour}, \ref{figure_bumpcontour}, we plot the contours of $f$ at some selected time. In Figure \ref{figure_bumpspa}, we show the spatial average of $f$. Those results agree with the benchmarks in the literature.

\section{Concluding Remarks}
\label{sec:conclusion}

In this paper, we develop Eulerian explicit and implicit solvers that can conserve total energy of the VA system. In particular, energy-conserving operator splitting is used for the fully implicit schemes. Numerical results demonstrate the accuracy, conservation, and robustness of our methods for three benchmark test problems. Our next goal is to generalize the methods to the VM system.

\section*{Acknowledgments}
 YC is supported by grants NSF DMS-1217563, AFOSR FA9550-12-1-0343 and the startup fund from Michigan State University. AJC is supported by AFOSR grants FA9550-11-1-0281, FA9550-12-1-0343 and FA9550-12-1-0455,  NSF grant DMS-1115709
and MSU foundation SPG grant RG100059.  We gratefully acknowledge the support from Michigan Center for Industrial and Applied Mathematics.

\end{document}